\documentclass[a4paper,final,3p,reqno,authoryear]{elsarticle}

\usepackage{hyperref}

\usepackage{array}

\usepackage{graphicx}

\usepackage{amsmath}
\usepackage{mathptmx}
\usepackage{color}
  \usepackage{bbm}
 \usepackage{amssymb}
  \usepackage{natbib}
\newcounter{Mtheorem}
\newenvironment{Mtheorem}[1][]{\refstepcounter{Mtheorem}\par\noindent\textbf{~Theorem~\theMtheorem. #1}\begin{itshape}}% begin code
{\end{itshape} \par}
\newcounter{trans}
{ \end{itshape} \par}
\newcounter{Mcorollary}
\newenvironment{Mcorollary}[1][]{\refstepcounter{Mcorollary}\par\noindent\textbf{~Corollary~\theMcorollary. #1}\begin{itshape}}% begin code
{\end{itshape} \par}
\newcounter{Mlemma}
\newenvironment{Mlemma}[1][]{\refstepcounter{Mlemma}\par\noindent\textbf{~Lemma~\theMlemma. #1}\begin{itshape}}% begin code
{\end{itshape}\par}
\newcounter{Proof}
{\end{itshape}\par}
\newcounter{Mexample}
\newenvironment{Mexample}[1][]{\refstepcounter{Mexample}\par\noindent\textit{~Example~\theMexample. #1}\begin{itshape}}% begin code
{  \end{itshape}\par}
\newcounter{Mdefinition}
\newenvironment{Mdefinition}[1][]{\refstepcounter{Mdefinition}\par\noindent\textbf{~Definition~\theMdefinition. #1}\begin{itshape}}% begin code
{ \end{itshape}\par}
\newcounter{Mproposition}
{ \end{itshape}\par}

\begin{document}

\begin{frontmatter}

\title{Imprecise Dirichlet process  with application  to the hypothesis test on the probability that $X\leq Y$}
% \tnotetext[mytitlenote]{Fully documented templates are available in the elsarticle package on \href{http://www.ctan.org/tex-archive/macros/latex/contrib/elsarticle}{CTAN}.}

\author[mymainaddress1]{Alessio Benavoli}
\author[mymainaddress1]{ Francesca Mangili}
\author[mymainaddress2]{Fabrizio Ruggeri}
\author[mymainaddress1]{Marco Zaffalon}

\address[mymainaddress1]{ IPG IDSIA, Manno, Switzerland}
\address[mymainaddress2]{ CNR IMATI, Milano, Italy}

\begin{abstract}
The Dirichlet process (DP) is one of the most popular Bayesian nonparametric models. 
An open problem with the DP is how to choose its infinite-dimensional parameter (base measure) in case of lack of prior information. 
In this work we present the Imprecise DP (IDP)---a prior near-ignorance DP-based model that does not require any choice of this probability measure.
It consists of a class of DPs obtained by letting the normalized base measure of the DP vary in the set of all probability measures. We discuss the tight connections of this approach with Bayesian robustness and in particular prior near-ignorance modeling via sets of probabilities.
We use this model to perform a Bayesian hypothesis test on the probability $P(X\leq Y)$. We study the theoretical properties of the IDP test (e.g., asymptotic consistency), and compare it
with the frequentist Mann-Whitney-Wilcoxon rank test that is commonly employed as a test on $P(X\leq Y)$. In particular we will show that our method is more robust, in the sense that it is able to isolate instances in which the aforementioned  test is virtually guessing at random.
\end{abstract}


\end{frontmatter}

\section{Introduction}
\label{sec:1}
The Dirichlet process (DP) is one of the most popular Bayesian nonparametric models. 
It was introduced by \citet{Ferguson1973} as a prior over probability distributions.
In his seminal paper, Ferguson showed that the DP leads to tractable posterior inferences
and can be used for Bayesian analysis of several nonparametric models, such as the estimation of a distribution function, of a mean, of quantiles, of a variance, etc. He also considered the estimation of  $P (X \leq Y )$ assigning independent Dirichlet process priors to the distribution functions of $X$ and $Y$. 
The Mann-Whitney statistic naturally arises in this case. 
\citet{SusarlaRyzin1976} and \citet{Blum1977} extended the results of Ferguson on
estimation of the distribution function in case of right censored data obtaining a Bayesian
version of the Kaplan-Meier estimator.
\citet{dalai1983nonparametric} considered the problem of
estimating  a measure of dependence for a
bivariate distribution. The Bayes estimate is computed using a two-dimension Dirichlet prior 
and Kendall's tau is seen to appear naturally.
A review of other similar applications of the DP can be found in \citep{Phadia2013}.

The beauty of the DP is that most of these results are in closed form and that it
provides a Bayesian justification of  the classic nonparametric estimators. In spite of all these nice properties and of the promising initial outcomes, it is a matter of fact that such a research  has not resulted in the development of  DP-based Bayesian nonparametric procedures for hypothesis testing.
For instance, the most used statistical packages for DP-based Bayesian nonparametric modeling, ``DPpackage'' \citep{Dppackage} and ``Bayesm'' \citep{rossi2010bayesm}, include procedures for density estimation, clustering and regression, but do not include any Bayesian version of the Wilcoxon rank sum, Wilcoxon sign test or other classic nonparametric tests.

It is arguable that this absence may be related to the unsettled question of how to choose the prior ``parameters'' of the DP in case of lack of prior information. 
It is well known in fact that a DP is completely characterized by its prior ``parameters'': the prior strength (or precision), which is a positive scalar number, and the normalized base measure. How should we  choose these prior ``parameters'' in case of lack of prior information? 
The only solution to this problem that has been proposed so far, first by \citet{Ferguson1973} and then by \citet{Rubin1981} under the name of Bayesian Bootstrap (BB), is the limiting DP obtained when the prior strength goes to zero. But the BB model has faced quite some
controversy, since it is not actually noninformative and moreover it assigns zero posterior probability to any set that does not include the observations. We will discuss these two points with more details in Section \ref{sec:PI}. 

In this paper we present an alternative viewpoint to the problem of choosing the prior base measure of the DP in case of lack of prior information that overcomes the above drawbacks of the BB. The model we present generalizes to nonparametrics  earlier ideas developed
in  Bayesian parametric robustness, see \citet{berger1994overview} and \citet{Ruggeri2000} for a review.
Here lack of prior information is expressed in terms of a family $\mathcal{T}$ consisting of all prior probability measures that are compatible with the available prior information. Inferences are then  carried out by considering the whole family $\mathcal{T}$.  In case almost no prior information is available on the parameters of interest, $\mathcal{T}$ should
be as large as possible in order to describe this state of prior ignorance.
The natural candidate for $\mathcal{T}$ to represent complete ignorance is the set of all probability measures.
However, it turns out that the posterior inferences obtained from this set are \textit{vacuous} \citep[Sec.~7.3.7]{walley1991}, i.e., 
 the posterior set   coincides with the prior set. This means that that there is no learning from data. Therefore, the vacuous prior model is not a statistically useful way to model our prior ignorance.
There is then a compromise to be made. \citet{pericchi_walley_91} and \citet{walley1991} suggest, as an alternative, the use of an almost vacuous model which they call ``near-ignorance'' or ``imprecise'' model. This is a model that behaves a priori as a vacuous model for some basic inferences (e.g., prior mean, prior credible regions)  but it always provides non-vacuous posterior inferences. 

While Bayesian robust models have already been extended to the nonparametric setting \citep{ruggeri1994}, 
that has not been the case for near-ignorance models.\footnote{A nonparametric model that uses lower and upper bounds for probabilities to quantify uncertainty
has been proposed by \citet{Augustin2004251,Coolen2009217}. This model is a purely predictive model, based on post-data assumptions, and,
thus, it cannot be used straightforwardly  to perform hypothesis tests.}
The main aim of this  paper is to derive a prior near-ignorance DP, called Imprecise DP (IDP).
This is the class $\mathcal{T}$ of all DPs obtained by fixing the prior strength of the DP and letting the normalized base measure vary in the set of all probability measures. We will show that the IDP  behaves a priori as a vacuous model for all predictive inferences.
This, together with the fact that it is a nonparametric model, allows us to start a statistical analysis with very weak assumptions about the problem of interest. However, contrarily to a full vacuous model,  we will show that the IDP can learn from data.

Moreover, we will employ the IDP to develop a new Bayesian nonparametric hypothesis test on the probability that $X\leq  Y$; we will call this test IDP rank-sum test, due to its similarity with the Mann-Whitney-Wilcoxon (MWW) rank-sum test.  
This hypothesis test is widely applied; for instance, if $X$ and $Y$ are health status measures in a clinical trial, $P(X\leq  Y )$ is, roughly speaking, the probability that the treatment represented by $Y$ is better (not worse) than the treatment represented by $X$. 
A Bayesian nonparametric near-ignorance model presents several advantages with respect to a traditional approach to hypothesis testing.
First of all, the Bayesian approach allows us to formulate the  hypothesis test as a decision problem. This means that we can  verify the evidence in favor of the null hypothesis and not only rejecting it, as well as take decisions that minimize the expected loss.
Second, because of the nonparametric near-ignorance prior, the IDP rank-sum test allows us to start the hypothesis test 
with very weak prior assumptions, much in the direction of letting data speak for themselves.
From a computational point of view, we will show that posterior inferences from the IDP can be derived by computing lower and upper bounds of expectations w.r.t.\
the class of DPs $\mathcal{T}$ and that,  for some inference, these lower and upper bounds can be computed in closed-form (e.g., mean and variance of $P(X\leq  Y )$).
When  no closed form expression exists,  these bounds can be computed by a simple Monte Carlo sampling from two  Dirichlet distributions. This means that we do not need to use stick breaking or other sampling approaches specific for DP. This computational advantage comes for free from prior near-ignorance.

In our view, the IDP rank-sum test appears to be a natural way to complete the work of \citet{Ferguson1973}, who first showed the connection between the expectation of 
$P(X\leq  Y )$ w.r.t.\ the DP and the Mann-Whitney statistic: it develops a Bayesian nonparametric near-ignorance-prior test for
the probability that $X\leq  Y $, which is computationally efficient and that, also for this reason, provides an effective practical alternative to the MWW  test.

Note that, although the IDP test shares several similarities with a standard Bayesian approach, at the same time it embodies a significant change of paradigm when it comes to take decisions.
In fact the IDP rank-sum test has the advantage of producing an \textit{indeterminate} outcome when the decision is \textit{prior-dependent}. In other words, the IDP test suspends the judgment (which can be translated as ``I do not know whether Y is  better than X'') when the option that minimizes the expected loss changes depending on the DP base measure we focus on. Therefore, the IDP-based test is \textit{robust} in the sense that it provides a determinate decision only when all the DPs, in the class the IDP represents, agree on the same decision.  
We will show that the number of indeterminate instances decreases at the accumulation of evidence and thus that the IDP-based test is always asymptotically consistent for $P(X\leq  Y )$.  This is not always true for the MWW test, even though the MWW test is commonly employed as a test about $P(X\leq  Y )$.

Finally, we will compare our IDP test with the MWW test and the DP-based test obtained as the prior strength goes to zero (called BB-DP test).
We empirically show on several different case studies that when the IDP test is indeterminate, 
the MWW and BB-IDP tests are virtually behaving as random guessers.
For a sample size of $20$ observations, the percentage of these instances can arrive up to almost $20\%$.
% This means that once every five decisions of MWW or BB-IDP tests are random guesses.
We regard this surprising result as an important finding, with practical consequences in hypothesis testing. Assume that we are trying to compare the effects of two  medical treatments (``Y  is better than X'') and that, given the available data, the IDP test is indeterminate. In such a situation the MWW test (or the BB-IDP test) always issues a determinate response (for instance, ``I can tell that Y  is better than X''), but it turns out that its response is completely random, like if we were tossing  a coin. In these cases by using MWW we would choose treatment $Y$, but this decision would be based on a random guess. In fact in these instances the MWW test could return the other hypothesis (``it is not true that Y  is better than X'') with equal probability.
On the other side, the IDP test acknowledges the impossibility of making a decision in these cases. Thus, by saying ``I do not know'', the IDP test  provides a richer information to the analyst. The analyst could for instance use this information to collect more data. (Please note that R and Matlab codes of our IDP rank-sum test are freely available at\\ \url{http://ipg.idsia.ch/software/IDP.php}.)

\section{Dirichlet process}
\label{sec:DP}
The Dirichlet process was developed by  \citet{Ferguson1973} as a probability distribution on the space of probability distributions. 
Let  $\mathbb{X}$ be a standard Borel space with Borel $\sigma$-field $\mathcal{B}_\mathbb{X}$  and $\mathbb{P}$ be the space
of probability measures on  $(\mathbb{X},\mathcal{B}_{\mathbb{X}})$ equipped with the weak topology and the
corresponding Borel $\sigma$-field $\mathcal{B}_\mathbb{P}$. Let $\mathbb{M}$ 
be the class of all probability measures on $(\mathbb{P},\mathcal{B}_\mathbb{P})$. We call the elements $\mu \in \mathbb{M}$
 nonparametric priors.

An element of  $ \mathbb{M}$ is called a Dirichlet process distribution $\mathcal{D}(\alpha)$ with base measure $\alpha$ if for every finite measurable
partition $B_1,\dots,B_m$ of $\mathbb{X}$, the vector $(P(B_1),\dots,P(B_m))$ has a Dirichlet distribution with parameters $(\alpha(B_1),\dots,\alpha(B_m))$, where $\alpha(\cdot)$ is a finite positive Borel measure on $\mathbb{X}$.
 Consider the partition $B_1=A$ and $B_2=A^c=\mathbb{X}\backslash A$ for some measurable set $A \in \mathbb{X}$, then if $P \sim \mathcal{D}(\alpha)$ from the definition of the DP we have that
$(P(A),P(A^c))\sim Dir(\alpha(A),\alpha(\mathbb{X})-\alpha(A) )$, which is a Beta distribution.
From the moments of the Beta distribution, we can thus derive that:
\begin{equation}
\label{eq:priormom}
\begin{array} {rl}
\mathcal{E}[P(A)]=\dfrac{\alpha(A)}{\alpha(\mathbb{X})},  &\mathcal{E}[(P(A)-\mathcal{E}[P(A)])^2]=\dfrac{\alpha(A)(\alpha(\mathbb{X})-\alpha(A))}{(\alpha(\mathbb{X})^2(\alpha(\mathbb{X})+1))},
\end{array}
\end{equation}
where we have used the calligraphic letter $\mathcal{E}$ to denote  expectation w.r.t.\ the Dirichlet process.
This shows that the normalized measure $\alpha(\cdot)/\alpha(\mathbb{X})$ of the DP reflects the prior expectation of $P$, while the scaling parameter $\alpha(\mathbb{X})$ controls how much $P$ is allowed to deviate from its mean $\alpha(\cdot)/\alpha(\mathbb{X})$.
Let $s = \alpha(\mathbb{X})$ stand for the total mass of $\alpha(\cdot)$  and $\alpha^*(\cdot)=\alpha(\cdot)/s$ stand for the
probability measure obtained by normalizing $\alpha(\cdot)$. If $P \sim \mathcal{D}(\alpha)$, we shall also
describe this by saying $P \sim Dp(s,\alpha^*)$ or, if $\mathbb{X} = \mathbb{R}$, $P \sim Dp(s, G_0)$, where $G_0$ stands
for the cumulative distribution function of $\alpha^*$.

Let $P\sim Dp(s,\alpha^*)$ and $f$ be a real-valued bounded function defined on $(\mathbb{X},\mathcal{B})$. Then the expectation with respect to the Dirichlet process of $E[f]$ is
\begin{equation}
\label{eq:expf}
\mathcal{E}\big[E(f)\big]=\mathcal{E}\left[\int f dP\right]=\int f d\mathcal{E}[P] = \int f d\alpha^*.
\end{equation}
One of the most remarkable properties of the DP priors is that the posterior distribution of $P$ is again a DP.
Let $X_1,\dots,X_n$ be an independent and identically distributed sample from $P$ and $P \sim Dp(s,\alpha^*)$, then the posterior distribution of $P$ given the observations is  
\begin{equation}
\label{eq:DPPost}
P|X_1,\dots,X_n \sim Dp\left(s+n, \frac{s}{s+n} \alpha^*+ \frac{1}{s+n} \sum\limits_{i=1}^n \delta_{X_i}\right),
\end{equation}
where $\delta_{X_i}$ is an atomic probability measure centered at $X_i$. This means that the Dirichlet process satisfies a property of conjugacy, in the sense that the posterior for $P$ is again a Dirichlet process with updated unnormalized base measure $\alpha+ \sum_{i=1}^n \delta_{X_i}$.
From (\ref{eq:DPPost}) and (\ref{eq:priormom})--(\ref{eq:expf}), we can easily derive the posterior mean and variance of $P(A)$ and, respectively, posterior expectation of $f$.
% This means that the Dirichlet process satisfies a property of conjugacy, in the sense that the posterior for $P$ is again a Dirichlet process with updated base measure.
Hereafter we list some useful properties of the DP that will be used in the sequel (see \citet[Ch. 3]{ghosh2003bayesian}).
\begin{description}
 \item[(a)] In case $\mathbb{X} = \mathbb{R}$, since $P$ is completely defined by its cumulative distribution function $F$, a-priori we say $F \sim Dp(s, G_0)$ and a posteriori we can rewrite (\ref{eq:DPPost}) as follows:
%  in an equivalent form. Assume that $F \sim Dp(s, G_0)$,  the posterior distribution of $F$ given $X_1,\dots,X_n$ is:
\begin{equation}
\label{eq:DPCDF}
F|X_1,\dots,X_n \sim Dp\left(s+n, \frac{s}{s+n} G_0+ \frac{1}{s+n} \sum_{i=1}^n I_{[X_i,\infty)}\right),
\end{equation} 
where $I$ is the indicator function.
\item[(b)] Consider an element $\mu \in \mathbb{M}$ which puts all its mass at the probability measure $P=\delta_{x}$ for some $x \in \mathbb{X}$.
This can also be modeled as  $Dp(s,\delta_x)$ for each $s>0$.
\item[(c)] Assume that $P_1 \sim Dp(s_1,\alpha_1^*)$, $P_2 \sim Dp(s_2,\alpha_2^*)$, $(w_1,w_2)\sim Dir(s_1,s_2)$ and $P_1$, $P_2$, $(w_1,w_2)$ are independent, then \citep[Sec. 3.1.1]{ghosh2003bayesian}:
\begin{equation}
 w_1P_1+w_2P_2 \sim Dp\left(s_1+s_2, \frac{s_1}{s_1+s_2} \alpha_1^*+ \frac{s_2}{s_1+s_2} \alpha_2^*\right).
\end{equation} 
\item[(d)] Let $P_x$ have distribution $Dp(s+n, \frac{s}{s+n} \alpha^*+ \frac{1}{s+n} \sum_{i=1}^n \delta_{X_i})$. We can write
\begin{equation} \label{eq:mixing}
P_x=w_0 P+ \sum_{i=1}^n w_i \delta_{X_i},
\end{equation} 
where  $(w_0,w_1,\dots,w_n)\sim Dir(s,1,\dots,1)$ and $P \sim Dp(s,\alpha^*)$ (it follows by (b)-(c)).
\end{description}

\section{Prior ignorance}
\label{sec:PI}
How should we choose the prior parameters $(s,\alpha^*)$ of the DP, in particular the infinite-dimensional  $\alpha^*$, in case of lack of prior information? 
To address this issue, the only prior that has been proposed so far is the limiting DP obtained for $s \rightarrow 0$, which has been introduced under the name of Bayesian Bootstrap (BB) by \citet{Rubin1981}; in fact it can be proven that the BB is asymptotically equivalent (see \citet{Lo1987} and \citet{Weng1989}) to the frequentist bootstrap introduced by  \citet{Efron1979}.

The BB has been criticized on diverse grounds. 
From an a-priori point of view, the main criticism is that taking $s \rightarrow 0$ is far from leading to  a noninformative prior.
\citet{sethuraman1981} have shown that for $s \rightarrow 0$ a measure sampled
from the DP is a degenerated (atomic) measure centered on $X_0$, with $X_0$ distributed according to $\alpha^*$.
% if we consider a sequence of measures $\alpha_r$ such that, $s \rightarrow 0$ and the normalized measure $\alpha_r/s$ converges in the strong sense to the probability measure $\alpha^*$, then the sequence of Dirichlet processes with parameter $\alpha_r$, converges, for $r\rightarrow \infty$, to a random degenerated measure centered on $X_0$, with $X_0$ distributed according to $\alpha^*$. 
% This result depends on the fact that small values of $s$ increase the probability of ties among samples drawn from $P\sim Dp(s,\alpha^*)$ and, thus, as $s$ goes to $0$ the DP tends to concentrate on a single value.
As a further consequence, from an a-posteriori point of view, this choice for the prior gives zero probability to the event that a future observation is different from the previously observed data. 
 \citet{Rubin1981} reports the following extreme example. Consider the probability that $X>C$ where $C$ is a value larger than the largest observed value of $X$, i.e., $X_{(n)}$. The standard BB and bootstrap methods estimate such probability to be $0$ with zero variance, which is untenable if $X$ can assume different values from the $n$ previously observed.
 Rubin also remarks that one should expect a probability that $X$ is greater than or equal to $X_{(n)}$ of about $1/(n+1)$.
This shows that a Dirichlet prior with $s \rightarrow 0$ implies definitely a very strong (and not always reasonable) information about $P$, and hence it cannot be considered a noninformative prior.
On the other side, if we choose a DP prior with $s>0$, the inferences provided by this model will be sensitive to the choice of the normalized measure $\alpha^*$. If, for example, we decide to assign a ``tail'' probability of $1/(n+1)$ to $X>X_{(n)}$, in agreement with Rubin's intuition, the inferences will be different if we assume that the tail probability is concentrated on  $X_{(n)}$ or if we assume that it is spread from $X_{(n)}$ to a very large value of $X$.

To answer to the initial question of this section, we propose the imprecise Dirichlet process (IDP). 
 The main characteristic of the IDP is that it does not require
 any choice of the normalized measure $\alpha^*$, it is a prior near-ignorance model  and solves the issues of the BB.
% method to choose  $(s,\alpha^*)$ in case of lack of prior information: a \textit{prior near-Ignorance DP}. 
% This model  does not require any choice of the normalized measure $\alpha^*$.
% IDP has been derived extending \citet{pericchi_walley_91,walley1991}'s parametric prior near-ignorance models to the DP.
Before introducing the IDP, it is worth to explain  what is a prior near-ignorance model with the example of a parametric model \cite[Sec. 5.3.1]{walley1991}.
% So before introducing this model in the DP setting, it is worth to review  what is a prior ignorance model in the parametric setting.
% Prior Ignorance models have been introduced by \citet{pericchi_walley_91,walley1991} in the context of Bayesian parametric.
% These are models that behave a priori as a vacuous models for some basic inferences (e.g., prior mean, prior credible regions)  but it always provides non-vacuous posterior inferences. 
% 
% A Prior Ignorance model simply consist in a set of probability measures
% 
% introduced a class of priors, called ``near-ignorance priors'' see also \citet{walley1991}.
% This is a model that behaves a priori as a vacuous model for some basic inferences (e.g., prior mean, prior credible regions)  but it always provides non-vacuous posterior inferences.  
\begin{Mexample}
Let $A$ be the event that a particular thumbtack lands pin-up at the next toss. Your information is that there
have been $m$ occurrences of pin up in $n$ previous tosses.
 Using a Bernoulli model, the likelihood
function generated by observing m successes in n trials is then proportional
to $\theta^n(1-\theta)^{n-m}$ where $\theta$ is the chance of pin-up.
To complete the model, we need to specify prior beliefs concerning the
unknown chance $\theta$.
We can use a conjugate Beta prior $p(\theta)=Be(\theta;\alpha,\beta)$,
 where $\alpha,\beta>0$ are the prior parameters of the Beta density.
 A-posteriori we have that $p(\theta|m,n)=Be(\theta;\alpha+m,\beta+n-m)$.
 Thus, the prior and posterior probabilities of $A$ are:
 $$
 P(A)=E[\theta]=t,~~~P(A|m,n)=E[\theta|m,n]=\frac{st+m}{s+n},
 $$
 where $s=\alpha+\beta$ is the prior strength and $t=\alpha/(\alpha+\beta)$ the prior mean.
 The problem is how to choose the parameters $s,t$ in case of lack of prior information.
  \citet[Ch. 5]{walley1991} proposes to use a   prior near-ignorance model.
A  near-ignorance prior model for this example is any set of priors which generates vacuous prior probabilities for the event
of interest $A$, i.e.,
$$
\underline{P}(A)=0,~~\overline{P}(A)=1,
$$
where $\underline{P},\overline{P}$ are lower and upper bounds for $P(A)$. These
vacuous probabilities reflect a complete absence of prior information
concerning $A$. For the Beta prior, since $P(A)=E[\theta]=t$,  the class of priors is simply:
 $$
 p(\theta)\in\left\{Be(\theta;st,s(1-t)):~~0<t<1\right\},
 $$
 for some fixed $s>0$, i.e., this is the set of priors obtained by considering all the Beta densities whose mean parameter $t$ is free to span  the interval  $(0,1)$. %This set of priors  has been called ``Imprecise Beta Model'' by \citet{walley1991.
 Posterior inferences from this model are derived by computing lower and upper posterior bounds;
in the case of event $A$ these bounds are:
 $$
\underline{P}(A|m,n)=\dfrac{m}{s+n},~~\overline{P}(A|m,n)=\dfrac{s+m}{s+n},
$$
where the lower is obtained for $t \rightarrow 0$ and the upper for $t \rightarrow 1$. 
% These both approach the observed relative frequency $m/n$ as the sample
% size $n$ increases. Behaviourally, this means that you are initially unwilling to
% bet on or against $A$ at any odds, because of the near-ignorance property,
% but, as observations are made, you become willing to bet at rates approaching
% the observed relative frequency.
% Note that the set of Beta densities satisfy prior ignorance for other inferences on $\theta$, e.g., non-central
% moments, credible intervals etc..
% Finally, observe  that $s$ is the only parameter that has to be selected by the modeller and determine the degree of ignorance (or robustness)
% of the model, i.e., the difference between the upper and lower bound
% $$
% \overline{P}(A|m,n)-\underline{P}(A|m,n)=\dfrac{s}{s+n}.
% $$
We point the reader to \citet{walley1996} for more details about this model and to \citet{Benavoli20121960}
for an extension of near-ignorance to one-parameter exponential families.
\end{Mexample}

\subsection{Imprecise Dirichlet process}
Before introducing the IDP, we give a formal definition of (nonparametric) prior ignorance
for predictive inferences.
Let $f$ be a  real-valued bounded function on $\mathbb{X}$, we call $E[f]=\int f dP$ a predictive inference about $X$; here $P$
is a probability measure on $(\mathbb{X},\mathcal{B}_{\mathbb{X}})$. Let $\mu \in \mathbb{M}$ be
a nonparametric prior on $P$ and  $\mathcal{E}_{\mu}[E(P)]$ the expectation of $E[f]$ w.r.t.\ $\mu$.
\begin{Mdefinition}
\label{def:PI}
A class of nonparametric priors  $\mathcal{T}\subset \mathbb{M}$ is called
a prior ignorance model for  predictive inferences about $X$,  if  for any  real-valued bounded function $f$ on $\mathbb{X}$ it satisfies:
% Let $f$ be a real-valued bounded function defined on $\mathbb{X}$; we say that
% a set of probability measures $\mathcal{G}$ is a prior ignorance model for $E[f]$ if it holds that:
 \begin{equation}
\label{eq:vacdef}
\underline{\mathcal{E}}[E(f)]=\inf\limits_{\mu \in \mathcal{T}} \mathcal{E}_{\mu}[E(f)] =\inf f, ~~ \overline{\mathcal{E}}[E(f)]=\sup\limits_{\mu \in \mathcal{T}} \mathcal{E}_{\mu}[E(f)]=\sup f,
\end{equation}
where  $\underline{\mathcal{E}}[E(f)]$ and $\overline{\mathcal{E}}[E(f)]$ denote respectively the lower and upper bound of $\mathcal{E}_{\mu}[E(P)]$ 
 calculated w.r.t.\ the class $\mathcal{T}$.
\end{Mdefinition}
From (\ref{eq:vacdef}) it can be observed that the range of $\mathcal{E}_{\mu}[E(f)]$  under the class $\mathcal{T}$ is
the same as the original range of $f$. 
In other words, by specifying the class  $\mathcal{T}$, we are not giving any information
on the value of the expectation of $f$. This means that  the class $\mathcal{T}$ behaves as a vacuous model. 
% All predictive inferences about a rv $X$ can be written as the expected value of a function $f(X)$. For example, for the probability $P(X\in A)$ for some $A \subseteq \mathbb{X}$ can be written as $E(f)$, with $f(x) = I_A(x)$. Given a nonparametric prior $\mu$, we can evaluate $E(f)$ by taking its expected value $\mathcal{E}_{\mu}[E(f)]$ w.r.t. $\mu$. A model of prior ignorance about $E(f)$ can only tell us that the expectation of $E(f)$ is included in the original range of $f$, i.e.,
% $$
% \inf{f} \leq \mathcal{E}_{\mu}[E(f)] \leq \sup{f}
% $$
% where we can replace infimum (supremum) with minimum (maximum) if the infimum (supremum) of $f$ is contained in the range of $f$. This would mean that the prior IDP model does not add any information about $E(f)$ and for this reason is a model of prior ignorance about the value of $E(f)$.
% Therefore, we give the following definition of prior ignorance in the nonparametric setting.
% \begin{Mdefinition}[Prior Ignorance.]
% We say that a set $\mathcal{T^*}$ of probability measures is a prior ignorance model for the predictive inference about $X$ denoted as $E(f)$, with $f$ a real-valued bounded function on $\mathbb{X}$, if it verifies 
% \begin{equation}
% \inf_{\mu \in \mathcal{T^*}} \mathcal{E}_{\mu}[E(f)]=\inf f, ~~ \sup_{\mu \in \mathcal{T^*}} \mathcal{E}_{\mu}[E(f)]=\sup f.
% \end{equation}
% \end{Mdefinition}
We are now ready to define the IDP.
\begin{Mdefinition}[IDP.]
We call  prior imprecise DP the following class of DPs:
 \begin{equation}
\label{eq:IDP}
\mathcal{T}=\left\{Dp(s,\alpha^*): ~\alpha^* \in \mathbb{P}\right\}.
\end{equation}
\end{Mdefinition}
The IDP is the class of DPs obtained for a fixed $s>0$ and by letting the normalized measure $\alpha^*$ to vary in the set of all 
probability measures $\mathbb{P}$ on  $(\mathbb{X}, \mathcal{B}_{\mathbb{X}})$.
\begin{Mtheorem}
\label{th:PI}
The IDP is a model of prior ignorance for all predictive inferences about $X$, i.e., for any  real-valued bounded function $f$ on $\mathbb{X}$ it satisfies:
% Let $f$ be a real-valued bounded function defined on $\mathbb{X}$; we say that
% a set of probability measures $\mathcal{G}$ is a prior ignorance model for $E[f]$ if it holds that:
 \begin{equation}
\label{eq:vac1}
\underline{\mathcal{E}}[E(f)]=\inf f, ~~ \overline{\mathcal{E}}[E(f)]=\sup f,
\end{equation}
where  $\underline{\mathcal{E}}[E(f)]$ and $\overline{\mathcal{E}}[E(f)]$ denote respectively the lower and upper bound of ${\mathcal{E}}[E(f)]$ 
defined in (\ref{eq:expf}) calculated w.r.t.\ the class of DPs (\ref{eq:IDP}).

% Let $f$ be a real-valued bounded function defined on $\mathbb{X}$; we say that
% a set of probability measures $\mathcal{G}$ is a prior ignorance model for $E[f]$ if it holds that:
 %\begin{equation}
%\label{eq:vac1}
%\underline{\mathcal{E}}[E(f)]=\inf f, ~~ \overline{\mathcal{E}}[E(f)]=\sup f,
%\end{equation}
%where  $\underline{\mathcal{E}}[E(f)]$ and $\overline{\mathcal{E}}[E(f)]$ denote respectively the lower and upper bound of ${\mathcal{E}}[E(f)]$ 
%defined in (\ref{eq:expf}) calculated w.r.t.\ the class $\mathcal{T}$ of DPs (\ref{eq:IDP}).

% on $\mathbb{X}$, where  
% $$
% \underline{\mathcal{E}}[E(f)]=\inf_{~\alpha^* \in \mathbb{P}} \mathcal{E}[E(f)|G], ~~~~\overline{\mathcal{E}}[E(f)]=\sup_{G \in \mathcal{G}} \mathcal{E}[E(f)|G],
% $$ 
% denote the lower and, respectively, upper  bound of  $\mathcal{E}[E(f)|G]$, which is the expectation of $E(f)$ with respect to some random probability measure $G$ in $\mathcal{G}$.
\end{Mtheorem}
The proofs of this and the next theorems are in the Appendix.  To show that the IDP is a model of prior ignorance, consider for instance the indicator function $f=I_{A}$ for some $A \subseteq \mathbb{X}$. Since $E[I_A]=P(A)$, from (\ref{eq:expf}) we have that $\mathcal{E}[P(A)]= \int I_A d\alpha^*$.
Then if we choose  $\alpha^*=\delta_{x_l}$ with $x_l\notin A$  and, respectively,  $\alpha^*=\delta_{x_u}$
with $x_u\in A$:
% Equation (\ref{eq:vac1}) becomes
 \begin{equation}
\label{eq:vac1ind}
\underline{\mathcal{E}}[P(A)]= \int I_A d\delta_{x_l}=\min I_A=0,~~ \overline{\mathcal{E}}[P(A)]= \int I_A  d\delta_{x_u}=\max I_A= 1,
\end{equation}
where $\underline{\mathcal{E}}[P(A)]$ and $\overline{\mathcal{E}}[P(A)]$ are the lower and upper
bounds for $\mathcal{E}\big[P(A)\big]$. This is a condition of prior ignorance for $P(A)$,
since we are saying that the only information about $P(A)$ is that $0\leq P(A) \leq 1$.
%  Note that in (\ref{eq:vac1ind}) we have replaced the infimum (supremum) of $f$ with minimum (maximum) because the infimum is contained in the range of $f$.
The lower and upper bounds are obtained from the degenerate DPs $Dp(s,\delta_{x_l})$ and $Dp(s,\delta_{x_u})$, which belong to the class (\ref{eq:IDP}).
Note that, although the lower and upper bounds are obtained by degenerate DPs, to obtain these bounds
we are considering all possible $Dp(s,\alpha^*)$ with $\alpha^* \in \mathbb{P}$ (even the ones with continuous probability measures $\alpha^*$).
\begin{Mtheorem}[Posterior inference.]
\label{th:DPpost}
Let $X_1,\dots,X_n$ be i.i.d.\ samples from $P$ and $P \sim Dp(s,\alpha^*)$. Then for any  real-valued bounded function $f$ on $\mathbb{X}$, 
the lower and upper bounds of ${\mathcal{E}}[E(f)|X_1,\dots,X_n]$ under the IDP model in (\ref{eq:IDP}) are:
% the posterior distribution of $P$ given the observations is 
% For any  real-valued bounded function $f$ on $\mathbb{X}$ it satisfies:
\begin{equation}
\label{eq:postexp}
\begin{aligned}
\underline{\mathcal{E}}\big[E(f)|X_1,\dots,X_n\big]&= \dfrac{s}{s+n}\inf f+\dfrac{n}{s+n}S_n(f),\\
\overline{\mathcal{E}}\big[E(f)|X_1,\dots,X_n\big]&=  \dfrac{s}{s+n}\sup f+\dfrac{n}{s+n}S_n(f),
\end{aligned}
\end{equation} 
where $S_n(f)=\frac{\sum_{i=1}^n f(X_i)}{n}$.
% and $\underline{\mathcal{E}}[E(f)|X_1,\dots,X_n],\overline{\mathcal{E}}[E(f)|X_1,\dots,X_n]$ denote the lower and upper bounds 
% of ${\mathcal{E}}[E(f)|X_1,\dots,X_n]$ w.r.t.\ the class of DPs (\ref{eq:IDP}).
\end{Mtheorem}
A-posteriori the IDP does not satisfy anymore the prior ignorance property (\ref{eq:vac1}). This means that learning from data takes place under the IDP. 
In fact let $S(f)$ be equal to $\lim_{n\rightarrow \infty} S_n(f)$, a-posteriori  for $n\rightarrow \infty$ we have that:
\begin{equation}
\begin{aligned}
\underline{\mathcal{E}}\big[E(f)|X_1,\dots,X_n\big],\overline{\mathcal{E}}\big[E(f)|X_1,\dots,X_n\big] &\rightarrow S(f),
\end{aligned}
\end{equation} 
i.e., the lower and upper bounds of the posterior expectations converge to $S(f)$, which only depends on data. 
In other words, the effect of prior ignorance vanishes asymptotically:
$$
\overline{\mathcal{E}}\big[E(f)|X_1,\dots,X_n\big]-\underline{\mathcal{E}}\big[E(f)|X_1,\dots,X_n\big]=\dfrac{s}{s+n} (\sup f-\inf f)\rightarrow 0, 
$$
for any finite $s$. To define the IDP, the modeler has only to choose $s$. 
This explains the meaning of the adjective \textit{near} in prior near-ignorance,  because the IDP requires by the modeller the elicitation of a parameter. However, this is a simple elicitation problem for a nonparametric prior, since we only have  to choose the value of a positive scalar (there are not infinitely dimensional parameters left in the IDP model). Section \ref{sec:MWW} gives some guidelines for the choice of this parameter.

% Fixed $s$, for any choice of $\alpha^*$ (continuous or discrete), we always have
% $$
% \underline{\mathcal{E}}\big[E(f)|X_1,\dots,X_n\big] \leq {\mathcal{E}}_{\alpha^*}\big[E(f)|X_1,\dots,X_n\big] \leq \overline{\mathcal{E}}\big[E(f)|X_1,\dots,X_n\big],
% $$
% which means that the IDP is robust, in the sense that encompasses all other posterior inferences obtained for any choice of $\alpha^*$.
% 
% 
% The modeller has only to choose $s$ (actually an upper bound for $s$ would be enough) to completely specify
% IDP. This is a simple elicitation problem for a nonparametric prior we have only to choose the value of a positive scalar
% (there are not infinitely dimensional parameters left). Section \ref{} gives some guidelines to the choice of this parameter.
% this model of prior ignorance satisfies the property of learning\footnote{The property of learning is satisfied if at least one of the posterior (upper and lower) estimates differs from the prior after a sufficiently large, but finite, number of data has been observed} for any bounded functions $f$. 
Observe that IDP solves the two main drawbacks of Bayesian Bootstrap. From the a-priori point of view, 
we have shown in (\ref{eq:vac1}) that the IDP is a model of prior ignorance for predictive inferences. Moreover, the prior distributions considered can assign a non-null probability to unobserved values of $X$. Then, considering Rubin's example about the probability that $X$ is greater than or equal to $X_{(n)}$, which is obtained as the expectation of $f=I_{[C,\infty)}$ with $C>X_{(n)}$, from (\ref{eq:postexp}) we have a-posteriori that $\underline{\mathcal{E}}\big[E(f)|data\big]=0$ and $\overline{\mathcal{E}}\big[E(f)|data\big]=\tfrac{s}{s+n}$. The upper expectation is greater than zero and, for $s=1$, it is equal to $1/(1+n)$. This result is obtained without specifying how the probability of $1/(1+n)$ is spread between the values $X>X_{(n)}$, and thus it is insensitive to the model specification of tail probabilities.
Note that the IDP reduces to the imprecise Dirichlet model proposed by  \citet{walley1996}, see also \cite{Bernard2005123,deCooman2009204}), when we limit ourselves to consider a finite measurable partition $B_1,\dots,B_m$ of $\mathbb{X}$. In this case, the set of priors $\{Dp(s,\alpha^*)$, $\alpha^* \in \mathbb{P}\}$, reduces to a set of Dirichlet distributions with parameters $(s\alpha^*(B_1),\dots, s\alpha^*(B_m))$.

\section{An application to hypothesis testing}
\label{sec:MWW}
Hypothesis testing is an important application of nonparametric statistics.   Recently there has been an increasing interest in the development of Bayesian nonparametric procedures for hypothesis testing.
For instance Bayesian nonparametric approaches to the two-sample problem have been proposed using Dirichlet process mixture models or (coupling-optional) Polya trees priors by \citet{borgwardt2009bayesian,holmes2009two,ma2011coupling,RePEc:eee:csdana}.
% In particular DPM or Polya tree are used to build a Bayes factor for testing whether two (or $k$ in case of \citet{RePEc:eee:csdana}) samples
% arise from the same distribution.
Although   prior near-ignorance  may also be useful in these models and in the two-sample problem, we do not follow this avenue in the present study.
Our focus is instead the hypothesis test $P(X\leq Y)\lesseqgtr P(X>Y)$ (equivalently  $P(X\leq Y)\lesseqgtr 0.5$), given independent random samples of sizes $n_1$ and $n_2$ from two populations. This problem arises, for example, if one wishes to compare the response $X$ of a population with respect to the response $Y$ of a different population in order to establish whether the two populations perform equally well or one population has generally ``better'' responses than the other.

The nonparametric test traditionally applied in such situations is the Mann-Whitney-Wilcoxon (MWW) rank-sum test.
The null hypothesis of the MWW rank-sum test is that the two populations are equal, that is, they come from the same distribution
$F_X(x) = F_Y (x)$. Let $X^{n_1} =\{X_1,\dots,X_{n_1}\}$ and $Y^{n_2} =\{Y_1,\dots,Y_{n_2}\}$ be two sequences of observations 
from the two populations. The MWW test is based on the fact that, if the two populations have the same distribution, the distribution of the linear rank statistic
\begin{equation}
 \label{eq:Umwm}
U=\sum\limits_{i=1}^{n_1}\sum\limits_{j=1}^{n_2} I_{[X_i,\infty)}(Y_j),
\end{equation} 
can be computed by considering all the possible random arrangements of the observations in $X^{n_1}$ and $Y^{n_2}$.
At the increase of $n_1$ and $n_2$, this distribution converges to a Normal distribution with mean and variance given by
\begin{equation}
 \label{eq:Umeanvar}
E\left[\frac{U}{n_1n_2}\right]=\dfrac{1}{2}, ~~Var\left[\frac{U}{n_1n_2}\right]=\dfrac{n_1+n_2}{12n_1n_2}.
\end{equation} 
It is worth to stress that $F_X(x) = F_Y (x)$ implies  $P(X\leq Y)=0.5$ (i.e., it is not true that $Y$ is better than $X$)
but not vice versa. Thus, the MWW test cannot be used in general as a test for $P(X\leq Y)$.
This limitation of the test is due to the choice of the  U statistic as estimator and the need of any frequentist method to specify the distribution of the  statistic under the null hypothesis. The null hypothesis $F_X(x) = F_Y (x)$ is thus selected  to be able to compute
the distribution of the statistic, although, in practice, one is interested in a much weaker hypothesis to test $P(X\leq Y)$ (see \citet{Fay2010} for a detailed discussion).
To overcome this issue of the MWW test, it is often common to assume a location-shift model, which states that the two populations can only differ in locations:
$F_Y(y) = F_X(y - \Delta)$. 
% The location-shift parameter $\Delta$ is also known as treatment effect: if population $X$ is taken as a control population and population $Y$ is the treatment population, then $\Delta$ is the expected effect due to the treatment. 
The goal is then to test the hypothesis that there is no treatment effect $\Delta = 0$ ($P(X\leq Y)=0.5$) versus the alternative $\Delta > 0$
($P(X\leq Y)>0.5$) or $\Delta < 0$ ($P(X\leq Y)<0.5$). 
% This assumption seems more realistic (although much stricter) than the previous one, and thus will be adopted in the numerical applications for fair comparisons of the frequentist and Bayesian tests. 
Under this assumption, the MWW test can be interpreted as a \citet{hodges1963} estimator. On the other side, the Bayesian approach provides the posterior distribution of $P(X\leq Y)$, which can be used to compute the probability of any hypothesis of interest. Therefore, we are not limited in the choice of the null hypothesis.
% 
% 
% However the MWW test has some drawback when applied for testing $P(X\leq Y)=0.5$. 
% Firstly, the most general and realistic alternative hypothesis to $F_X(x) = F_Y (x)$ is that $F_X(x) \neq F_Y (x)$, so that all possible distributions $F_X$ and $F_Y$ are included. Instead, when comparing the null hypothesis $F_X(x) = F_Y (x)$ against, for instance, the alternative hypothesis $P(X\leq Y)>0.5$, we should consider only the  set of all distributions $F_X$ and $F_Y$ except those that do not verify neither the null nor the alternative hypotheses (i.e., verify $P(X<Y)\leq 0.5$ but are not equal) (see \citet{Fay2010} for a detailed discussion). This limitation is due to the fact that the frequentist approach needs to specify the distribution of the U statistic under the null hypothesis, an thus the strong null hypothesis $F_X(x) = F_Y (x)$ is required, although, in practice, one is interested in a much weaker hypothesis. On the other side, the Bayesian approach provides the posterior distribution of $P(X\leq Y)$ based on which the probability of any hypothesis of interest can be computed. Therefore, we are not limited in the choice of the null hypothesis.

Another well-known drawback of the MWW test, which is common to all frequentist tests, is the controversial meaning of the p-value. The Bayesian approach to decision making allows basing the decisions on the value of the expected loss, whose practical meaning is much more intuitive. 
For example, the hypothesis test:
$$
\stackrel{\underbrace{P(X\leq Y)\leq  P(X > Y)}}{P(X\leq Y)\leq  0.5}
 \textit{ vs. } 
\stackrel{\underbrace{P(X\leq Y)>  P(X > Y)}}{P(X\leq Y) >  0.5}
$$
can be performed in a Bayesian way in two steps.
First we define a loss function
\begin{equation}
\label{eq:loss}
 L(P,a) =\left\{\begin{array}{ll}
  K_0 I_{\{P(X\leq Y)> 0.5\}} &\text{ if } a = 0,\\
  K_1 I_{\{P(X\leq Y)\leq  0.5\}} &\text{ if } a = 1.\\
 \end{array}\right.
\end{equation} 
The first row gives the loss we incur by taking the action $a=0$ (i.e., declaring that $P(X\leq Y)\leq  0.5$) when actually $P(X\leq Y) >  0.5$,
while the second row gives the loss we incur by taking the action $a=1$ (i.e., declaring that $P(X\leq Y)> 0.5$) when actually $P(X\leq Y)\leq 0.5$.
Second, we compute the expected value of this loss:
\begin{equation}
\label{eq:exploss}
\mathcal{E}\left[ L(P,a)\right]
 =\left\{\begin{array}{ll}
  K_0 \mathcal{P}\left[P(X\leq Y)> 0.5\right]&\text{ if } a=0,\\
  K_1 \mathcal{P}\left[P(X\leq Y)\leq   0.5\right]&\text{ if } a=1,\\
 \end{array}\right.
\end{equation}
where we have used the calligraphic letter $\mathcal{P}$ to denote the probability w.r.t. the DPs $F_X$ and $F_Y$.
Thus, we choose $a=1$ if
\begin{equation}
\label{eq:exploss_dec}
  K_1 \mathcal{P}\left[P(X\leq Y)\leq   0.5\right] \leq   K_0 \mathcal{P}\left[P(X\leq Y)>  0.5\right] \Rightarrow \mathcal{P}\left[P(X\leq Y)>  0.5\right]> \dfrac{K_1}{K_1+K_0},
\end{equation}
%Since $ \mathcal{P}\left[P(X\leq Y)\leq   \frac{1}{2}\right]=1-\mathcal{P}\left[P(X\leq Y)>  \frac{1}{2}\right]$, this  is equivalent to:
%\begin{equation}
%\label{eq:exploss_dec}
%\mathcal{P}\left[P(X\leq Y)>  \frac{1}{2}\right]> \dfrac{K_1}{K_1+K_0}.
%\end{equation}
or $a=0$ otherwise.
When the above inequality is satisfied, we can declare that $P(X\leq Y)>0.5$ with probability $\frac{K_1}{K_1+K_0}=1-\gamma$  (e.g., $1-\gamma=0.95$).
%Note that if $F_X$ and $F_Y$ are symmetric distributions, $P(X\leq Y)=0.5$ implies that they have the same median. Thus, the test $P(X\leq Y)>0.5$ is a valid test for the median under the location-shift model. On the other side, for different asymmetric distributions $F_X, F_Y$, it may be that $P(X\leq Y)=0.5$ but $F_X$and $F_Y$ have different medians. Thus, $P(X\leq Y)>1/2$ is not a test for the difference in median, but a test on whether ``Y is greater than X''.

Finally, based on the imprecise DP model developed in this paper, we can perform a Bayesian nonparametric test that, besides overcoming the limitation of the frequentist one, is based on extremely weak prior assumptions, and easy to elicit, since it requires only to choose the strength $s$ of the DP instead of its infinite-dimensional parameter $\alpha$. When using the IDP set of priors, we consider for $F_X$ and $F_Y$ all the possible DP priors with strength lower than or equal to $s$ (since all inferences obtained for $s'<s$ are encompassed by those obtained for $s$). All these priors give a posterior probability $\mathcal{P}\left[P(X\leq Y)> 0.5\right]$ included between the lower and upper bounds $\underline{\mathcal{P}}\left[P(X\leq Y)> 0.5\right]$ and $\overline{\mathcal{P}}\left[P(X\leq Y)>  0.5\right]$.
%We can employ the results of Theorem \ref{th:1} and Corollary \ref{co:1} to perform a hypothesis test
%based on the model of prior ignorance developed in this section.
%We thus select $a=1/2$ and, according to the decision (\ref{eq:exploss_dec}) for some $K_0,K_1$ , we verify if
Thus, according to the decision rule in (\ref{eq:exploss_dec}) for some $\gamma = \tfrac{K_0}{K_0+K_1}$, we verify if
$$
\underline{\mathcal{P}}\big[P(X\leq Y)> 0.5|X^{n_1},Y^{n_2}\big]>1-\gamma, ~~\overline{\mathcal{P}}\big[P(X\leq Y)> 0.5|X^{n_1},Y^{n_2}\big]>1-\gamma,
$$
and then:
\begin{enumerate}
 \item if both the inequalities are satisfied we can declare that $P(X\leq Y)$ is greater than $0.5$ with probability larger than $1-\gamma$;
 \item if only one of the inequality is satisfied (which has necessarily to be the one for the upper), we are in an indeterminate situation, i.e., we cannot decide;
\item if both are not satisfied, we can declare that the probability that $P(X\leq Y)$ is greater than $0.5$ is lower than the desired probability of $1-\gamma$.
\end{enumerate}
When our model of prior ignorance returns an indeterminate decision, it means that the evidence from the observations is not enough to declare either that the probability of the hypothesis being true is larger or smaller than the desired value $1-\gamma$; more measurements are necessary to take a decision.

\begin{figure}[h]
        \centering
        
         \begin{tabular}{cc}
                     \includegraphics[width=0.35\textwidth]{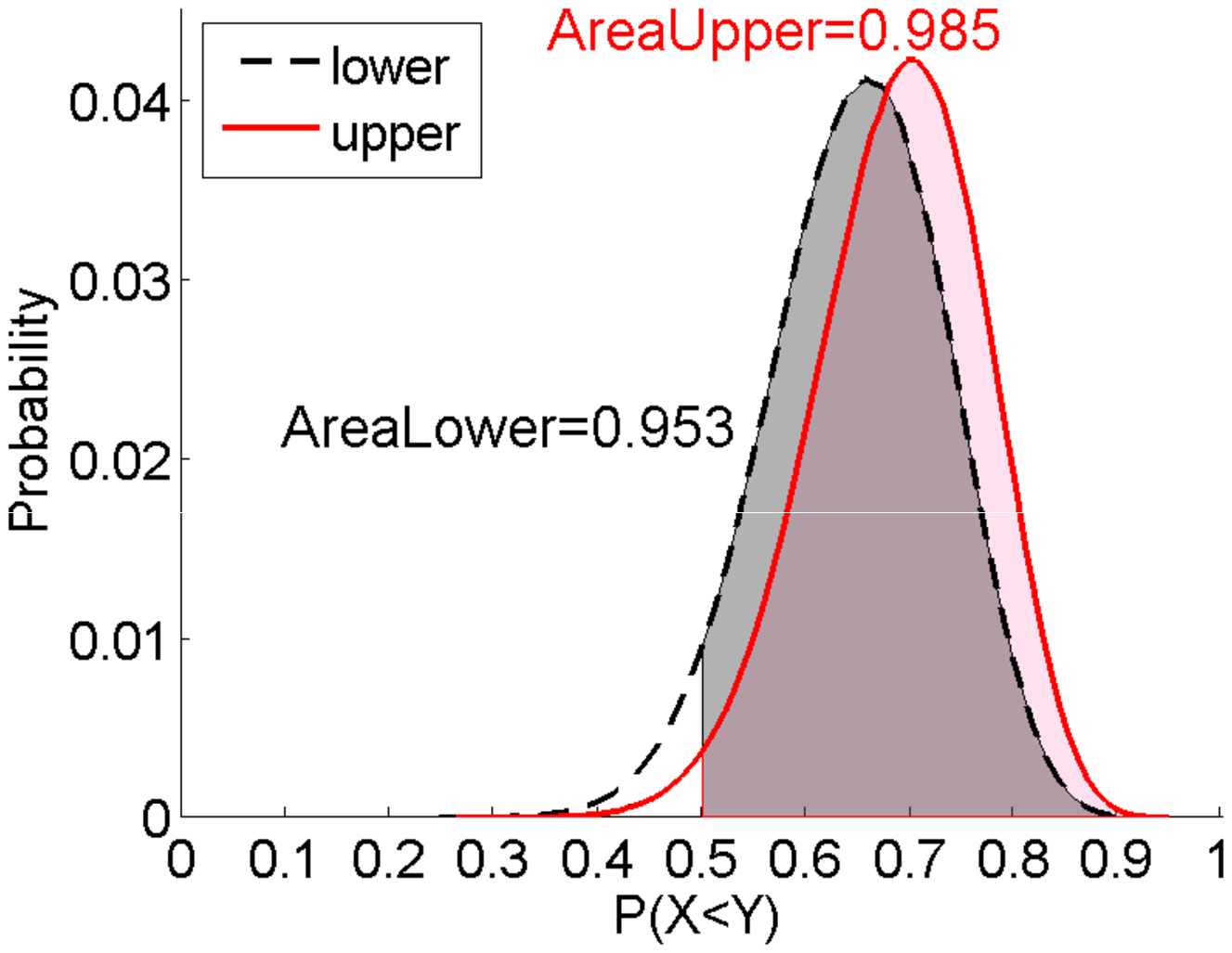} &  \includegraphics[width=0.35\textwidth]{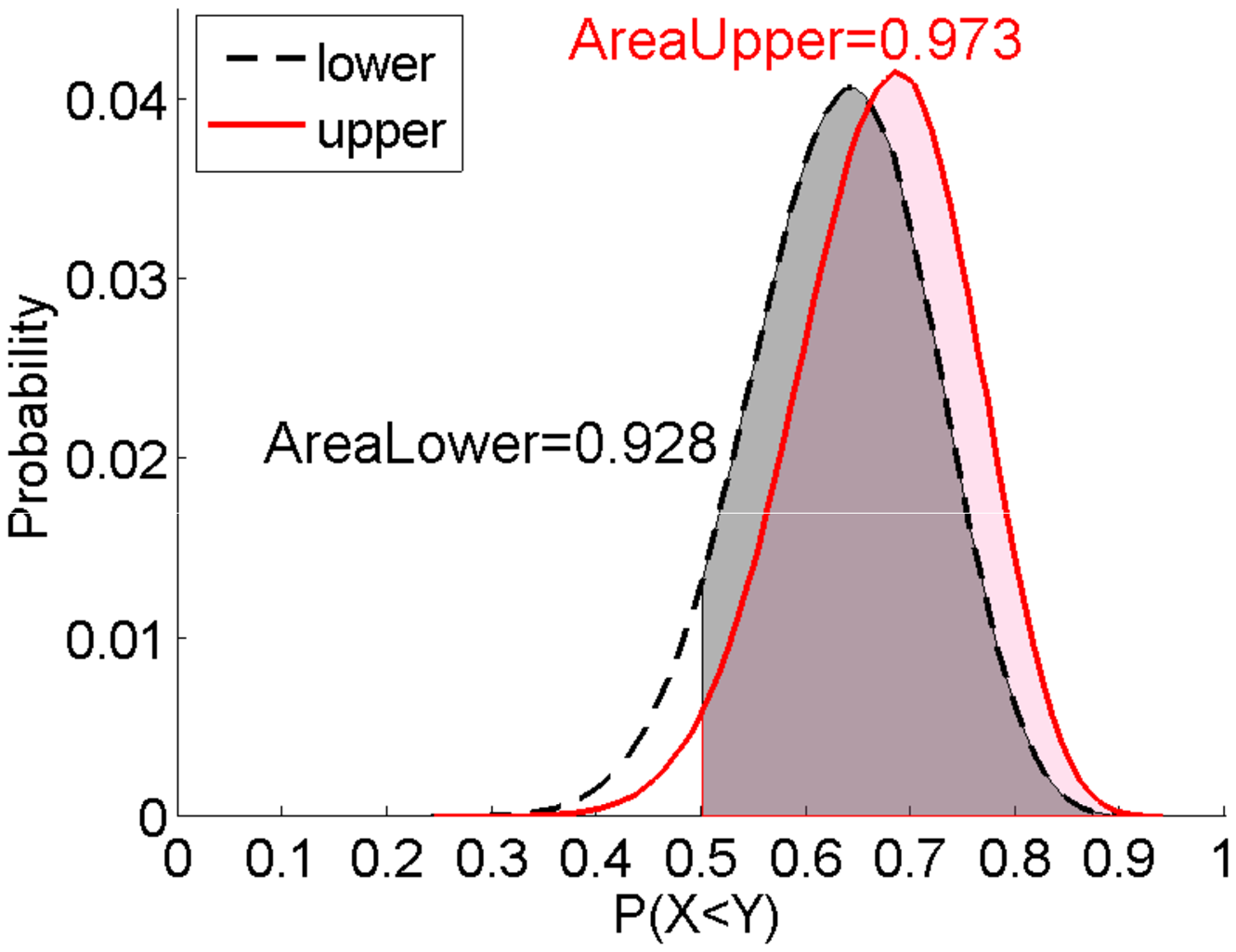}\\
                      (a) ``$Y$ is greater than $X$'' at $95\%$ &  (b) ``Indeterminate'' at $95\%$\\
                      \includegraphics[width=0.35\textwidth]{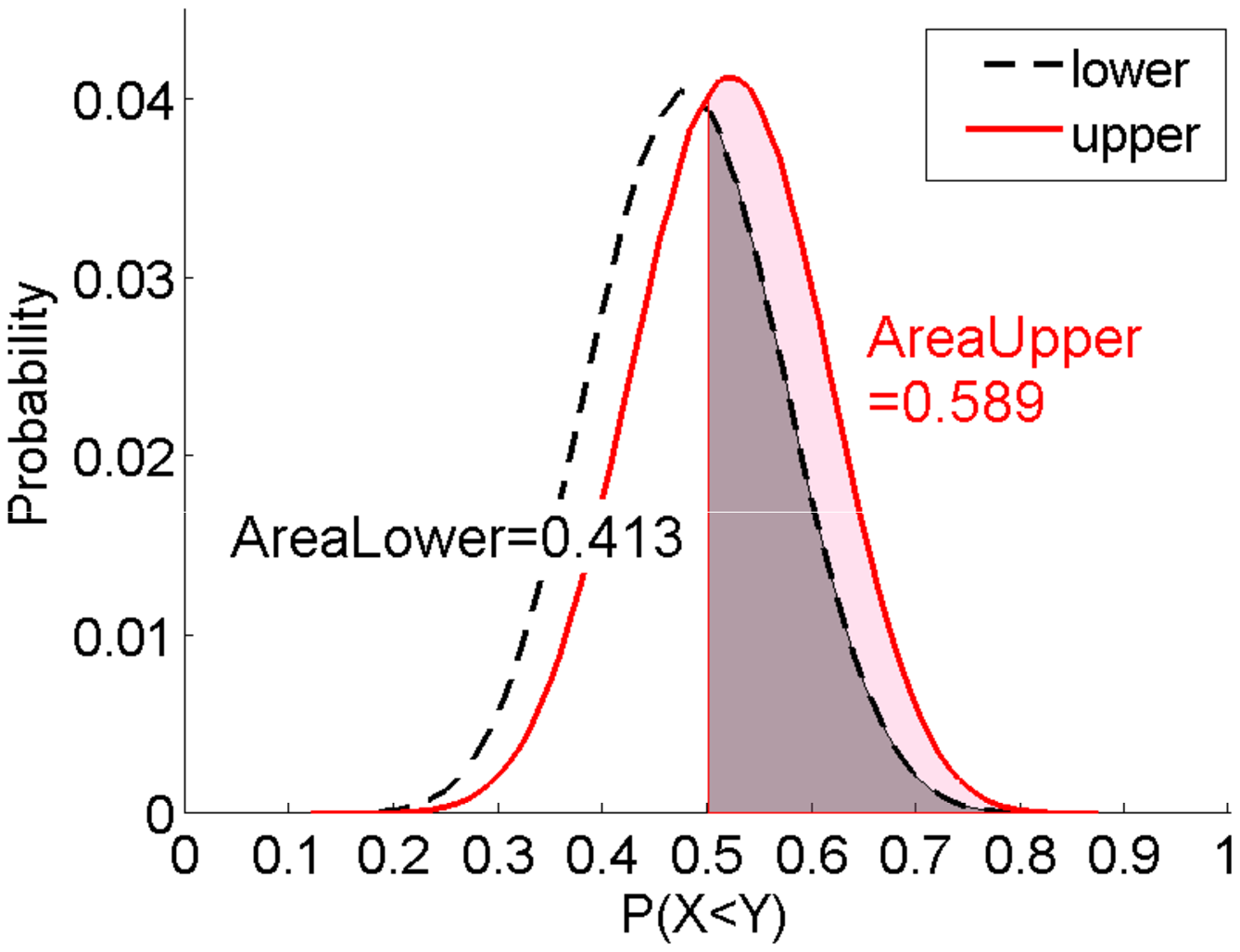} & \includegraphics[width=0.35\textwidth]{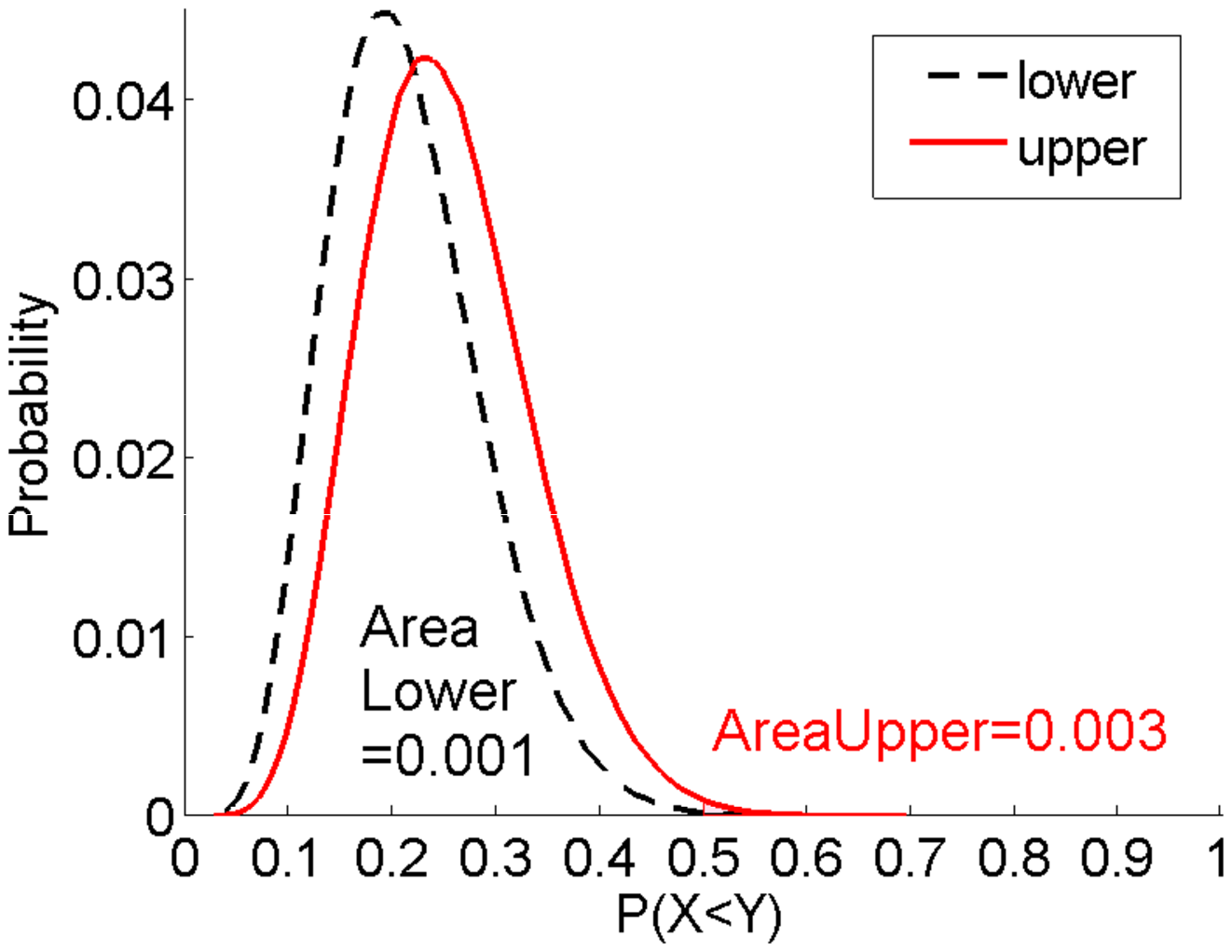}\\
                      (c) ``$Y$ is not greater than $X$'' at $95\%$ &  (d) ``$Y$ is not greater than $X$'' at $95\%$\\
         \end{tabular}
% 
%         \begin{subfigure}[b]{0.35\textwidth}
%                 \centering
%                 \includegraphics[width=\textwidth]{Test_determinate}
%                 \caption{T: ``$Y$ is greater than $X$'' at $95\%$}
%                 \label{fig:det}
%         \end{subfigure}%
%         ~ %add desired spacing between images, e. g. ~, \quad, \qquad etc.
%           %(or a blank line to force the subfigure onto a new line)
%         \begin{subfigure}[b]{0.35\textwidth}
%                 \centering
%                 \includegraphics[width=\textwidth]{Test_indeterminate}
%                 \caption{T: ``Indeterminate'' at $95\%$}
%                 \label{fig:ind}
%         \end{subfigure}
%         ~ %add desired spacing between images, e. g. ~, \quad, \qquad etc.
%           %(or a blank line to force the subfigure onto a new line)
%         \begin{subfigure}[b]{0.35\textwidth}
%                 \centering
%                 \includegraphics[width=\textwidth]{Test_determinateno0}
%                \caption{T: ``$Y$ is not greater than $X$'' at $95\%$}
%                 \label{fig:detn1}
%         \end{subfigure}
%         \begin{subfigure}[b]{0.35\textwidth}
%                 \centering
%                 \includegraphics[width=\textwidth]{Test_determinateno1}
%                \caption{T:``$Y$ is not greater than $X$'' at $95\%$}
%                 \label{fig:detn0}
%         \end{subfigure}
        \caption{Four possible results of the hypothesis test. The dark and light filled areas correspond respectively to the lower and upper probabilities of 
the event  $P(X\leq Y)>0.5$. The numerical values of these lower and upper probabilities are also reported in the figures. }\label{fig:1}
\end{figure}

The three cases are respectively depicted  in Figure  \ref{fig:1}.
Observe that the posterior  distributions  of $P(X\leq Y)$, from which the lower and upper probabilities above are derived, give us much more information than the simple result of the hypothesis test. In particular we can derive the posterior lower and upper probabilities of $P(X\leq Y)<0.5$.
For instance, from Figure \ref{fig:1} (d) we can see that $Y$ is not greater than $X$ at $95\%$, but it is evident that $X$ is  greater than $Y$ at $95\%$.
While in the case shown in Figure \ref{fig:1} (b), we can say neither that  $Y$ is  greater than $X$ nor that $X$ is  greater than $Y$. (To distinguish these two cases
it would be more appropriate to perform a ``two-sided'' hypothesis test.)

In the next section we prove that the IDP is a model of prior ignorance for $P(X\leq Y)$ and derive the posterior results which are necessary to evaluate $\mathcal{P}\left[P(X\leq Y)>  0.5\right]$ and perform the test. Note that, for the moment, we assume that there are not ties between $X$ and $Y$; we will discuss how to account for the presence of  ties in  Section \ref{sec:ties}.

%Ferguson \citet{Ferguson1973} has shown that the U statistic used in the MWW test appears naturally when computing the posterior expectation of $P(X<Y)$ given a DP prior; in particular, the posterior expectation of  $P(X<Y)$ converges to the U statistic as $s \rightarrow 0$. 

\subsection{IDP model for $P(X\leq Y)$}
\label{sec:MWWDP}
Let the samples $X^{n_1}$ and $Y^{n_2}$ be drawn, respectively, from $F_X$ and $F_Y$. As prior for $(F_X,F_Y)$, we assume that $F_X \sim Dp(s_1,G_1)$ and $F_Y \sim Dp(s_2,G_2)$. Hereafter, to simplify the presentation, we take $s_1 = s_2 = s$.
$F_X$ and $F_Y$ are assumed to be independent. 
The probability $P(X\leq Y)$ is given by $P(X\leq Y) =E[I_{[X,\infty)}(Y)] = \int F_X(y) dF_Y(y)$.
As derived by  \citet{Ferguson1973}, by the properties of the Dirichlet process, it follows that a-priori
$\mathcal{E}[P(X\leq Y)] = \int G_1(y) dG_2(y)$. It can be shown that the set of priors $\mathcal{T}$ in (\ref{eq:IDP}) satisfies the condition of prior ignorance also for $P(X\leq Y)$. In fact, since
$\mathcal{E}[P(X\leq Y)]=\int G_1(y)dG_2(y)$, if $G_i \in \mathbb{P}$, we have that
$$
\underline{\mathcal{E}}[P(X\leq Y)]=0, ~~\overline{\mathcal{E}}[P(X\leq Y)]=1,
$$
where the lower (upper) bound is obtained for  $dG_1 = \delta_{X_0}$ and $dG_2 = \delta_{Y_0}$ with $X_0>Y_0$ ($X_0<Y_0$).
Thus, prior ignorance about the mean of $P(X\leq Y)$ is satisfied.
Furthermore, let us consider the probability  of $P(X\leq Y)<0.5$ with respect to the Dirichlet process.  A-priori, for $dG_1 = \delta_{X_0}$ and $dG_2 = \delta_{Y_0}$ we have that
$$
\begin{array} {lll}
\text{if } X_0<Y_0, &\text{ then } \mathcal{P}[P(X\leq Y)=1]=1 &\text{ and thus } \underline{\mathcal{P}}[P(X\leq Y)\leq 0.5]=0\\
\text{if } X_0>Y_0, &\text{ then } \mathcal{P}[P(X\leq Y)=0]=1 &\text{ and thus } \overline{\mathcal{P}}[P(X\leq Y)\leq 0.5]=1. 
\end{array}
$$
A similar reasoning leads to $\underline{\mathcal{P}}[P(X\leq Y)>0.5]=0$, $\overline{\mathcal{P}}[P(X\leq Y)>0.5]=1$, 
thus, prior ignorance about the hypothesis $P(X\leq Y)>0.5$ is also satisfied.
Given the two sequences of measurements, a-posteriori one has:
$$
\mathcal{E}[P(X\leq Y)|X^{n_1},Y^{n_2}] = \int G_{n_1}(y) dG_{n_2}(y),
$$
with $G_{n_i}=\frac{s}{s+n_i} G_i+\frac{1}{s+n_i}\sum_{j=1}^{n_i} I_{[Z_j,\infty)}$, 
where $Z_j=X_j$ for $i=1$ and $Z_j=Y_j$ for $i=2$. It follows that:
 \begin{equation}
\label{eq:expmwmDp}
\begin{array}{rcl}
\mathcal{E}[P(X\leq Y)|X^{n_1},Y^{n_2}] &=& \frac{s}{s+n_1}\frac{s}{s+n_2}\int G_1(y) dG_2(y)+\frac{n_1}{s+n_1}\frac{s}{s+n_2}\frac{1}{n_1}\sum\limits_{j=1}^{n_1} (1-G_2(X_j^-))\\
 &+& \frac{s}{s+n_1}\frac{n_2}{s+n_2}\frac{1}{n_2}\sum\limits_{j=1}^{n_2} G_1(Y_j)+\frac{n_1}{s+n_1}\frac{n_2}{s+n_2}\frac{U}{n_1n_2},\\
\end{array}
\end{equation}
where $1-G_2(X^-)=\int I_{[X,\infty)}dG_2$. 
%The main issue of this approach is how  to select the prior distributions $G_1,G_2$ of the two Dirichlet processes for $X$ and $Y$ in case of lack of prior information. We have seen that in the Bayesian bootstrap (BB) ($s \rightarrow 0$) it holds that $G_{n_i}\rightarrow G^*_{n_i}=\tfrac{1}{n_i}\sum_{j=1}^{n_i} \delta_{Z_j}$, and thus the effect of the prior is eliminated from the posterior. Then the use of $G^*_{n_i}$ as a non-subjective (non-informative) posterior, is considered as one of the few concrete non-informative analysis in Bayesian nonparametrics \citet{ghosh2003bayesian}. However, as already discussed in Section \ref{sec:BB}, it is not completely correct to consider the BB prior as non-informative.
%\subsection{Modelling prior ignorance in the rank sum test}
%\label{sec:Modelling prior ignorance in the rank sum test}
%A way to model prior ignorance in the DP-based rank sum test is by following the approach discussed in Section \ref{sec:pign}.
%In fact, 
Then, the lower and upper posterior bounds of the posterior expectations of $P(X\leq Y)$ given the set of priors $\mathcal{T}$ are:
 \begin{equation}
\label{eq:ranktest}
\begin{array}{rcl}
\underline{\mathcal{E}}[P(X\leq Y)|X^{n_1},Y^{n_2}] &=&\frac{U}{(s+n_1)(s+n_2)},\\
\overline{\mathcal{E}}[P(X\leq Y)|X^{n_1},Y^{n_2}] &=& \frac{U}{(s+n_1)(s+n_2)}+\frac{s(s+n_1+n_2)}{(s+n_1)(s+n_2)}, \\
\end{array}
\end{equation}
obtained in correspondence of the extreme distributions $dG_1 \rightarrow \delta_{X_0}$, $dG_2 \rightarrow \delta_{Y_0}$,
with $X_0>\max(\{Y_0,\dots,Y_{n_1}\})$, $Y_0<\min(\{X_0,\dots,X_{n_2}\})$ (lower) and $X_0<\min(\{Y_0,\dots,Y_{n_1}\})$,  \\ $Y_0>\max(\{X_0,\dots,X_{n_2}\})$ (upper).
% This result provides lower and upper bounds for the expectation of $P(X\leq Y)$ derived by  \citet{Ferguson1973}.
% In this paper we complete the work of \citet{Ferguson1973} by develop
% For any other choice of $G_1,G_2$, it results that
% $$
% \underline{\mathcal{E}}[P(X\leq Y)|X^{n_1},Y^{n_2}]\leq {\mathcal{E}}[P(X\leq Y)|X^{n_1},Y^{n_2}] \leq \overline{\mathcal{E}}[P(X\leq Y)|X^{n_1},Y^{n_2}].
% $$
The posterior probability distribution of $P(X\leq Y)$ w.r.t. the Dirichlet process, which is used to perform the Bayesian test of the difference between the two populations, is, in general, computed numerically (Monte Carlo sampling) by using the stick-breaking construction of the Dirichlet process. We will show in the remaining part of this section that, in correspondence of the discrete priors that give the upper and lower bounds of the posterior distributions of $P(X\leq Y)$, a more efficient procedure can be devised. Consider the limiting posteriors that give the posterior lower and upper expectations in (\ref{eq:ranktest}):
\begin{equation}
\label{eq:basePost}
G_{n_i}(y)=\frac{s}{s+n_i} I_{[Z_{0},\infty)}+\frac{1}{s+n_i}\sum_{j=1}^{n_i} I_{[Z_{j},\infty)},
\end{equation}
where the lower is obtain with $Z_{0}=X_0>\max(\{Y_0,\dots,Y_{n_1}\})$ for $i=1$ and $Z_{0}=Y_0<\min(\{X_0,\dots,X_{n_2}\})$ for $i=2$, and the upper with $Z_{0}=X_0>\max(\{Y_0,\dots,Y_{n_1}\})$ for $i=1$, and $Z_{0}=Y_0<\min(\{X_0,\dots,X_{n_2}\})$ for $i=2$. 

% From the property (d) in Section \ref{sec:DP}, we can then derive the following result. 
\begin{Mlemma}\label{ml:1}
A sample $F_{n_i}$ from the Dirichlet process $Dp(s,G_{n_i})$ with base probability distribution $G_{n_i}$ as that defined in (\ref{eq:basePost}) is given by:
\begin{equation}
\label{eq:sample}
F_{n_i}=w_{i0} I_{[Z_{0},\infty)} + \sum\limits_{j=1}^{n_i} w_{ij} I_{[Z_{j},\infty)},
\end{equation}
where $(w_{i0},w_{i1},\dots,w_{in_i})\sim Dir(s,\stackrel{n_i}{\overbrace{1,\dots,1}})$.
\end{Mlemma}
Since $G_{n_i}$ is a discrete measure, Lemma \ref{ml:1} states that any distribution $F_{n_i}$ sampled from $DP(s+n_i,G_{n_i})$ has the  form (\ref{eq:sample}). Since the probability density functions relative to $F_{n_i}$, i.e., $w_{i0} \delta_{Z_{0}} + \sum\limits_{j=1}^{n_i} w_{ij} \delta_{Z_{j}}$, has a discrete support, we do not need stick-breaking. In fact, once we have selected the support $Z_0,Z_1,\dots,Z_{n_i}$, samples $F_{n_i}$ can simply be obtained by  sampling the weights from the Dirichlet distribution.
Using this fact, we can  derive lower and upper bounds for  ${\mathcal{P}}\big[P(X\leq Y)> c|X^{n_1},Y^{n_2}\big]$
as follows.

\begin{Mtheorem}
\label{th:1}
For any $c\in[0,1]$, it holds that 
\begin{equation}
\label{eq:lowmeandist}
\begin{array}{rcl}
&~& \underline{\mathcal{P}}\big[P(X\leq Y)> c|X^{n_1},Y^{n_2}\big]=P\left[g(w_{1\cdot},w_{2\cdot},X^{n_1}, Y^{n_2})>c\right],\\
\end{array}
\end{equation} 
with $$
g(w_{1\cdot},w_{2\cdot},X^{n_1}, Y^{n_2})=\sum\limits_{j=1}^{n_1}\sum\limits_{k=1}^{n_2} w_{1j}w_{2k} I_{(X_j,\infty)}(Y_k)
$$
% \begin{equation}
% \label{eq:lowmeandist}
% \begin{array}{rcl}
% &~& \underline{\mathcal{F}}\big[P(X\leq Y)\leq|X^{n_1},Y^{n_2}\big]:=\underline{\mathcal{P}}\big[P(X\leq Y)\leq |X^{n_1},Y^{n_2}\big]\vspace{2mm}\\
% &= &\displaystyle{\int \left[w_{10}I_{(X_0,\infty)}(y)+\sum_{j=1}^{n_1} w_{1j} I_{(X_j,\infty)}(y)\right]\left[w_{20}\delta_{Y_0}(y)+\sum_{j=1}^{n_2} w_{2j} \delta_{Y_j}(y)\right]dy}\\
% &=& \sum\limits_{j=1}^{n_1}\sum\limits_{k=1}^{n_2} w_{1j}w_{2k} I_{(X_j,\infty)}(Y_k),\\
% \end{array}
% \end{equation} 
where  $(w_{i0},w_{i1},\dots,w_{in_i})\sim Dir(s,\stackrel{n_i}{\overbrace{1,\dots,1}})$ for $i=1,2$ and $P$ is computed w.r.t.\ these Dirichlet distributions.
The mean and variance of $g(w_{1\cdot},w_{2\cdot},X^{n_1}, Y^{n_2})$ are:
\begin{equation}
\label{eq:lowmeanvar}
 \mu=E_W[W]^TAE_V[V], ~~~\sigma^2= trace[A^TE_W[WW^T]AE_V[VV^T]]-\mu^2,
\end{equation} 
where the expectations $E_W,E_V$ are taken w.r.t.\ the Dirichlet distributions, with  $W=[w_{11},\dots,w_{1n_1}]^T$, $V=[w_{21},\dots,w_{2n_2}]^T$, $E[WW^T]$ and $E[VV^T]$ are $n_i\times n_i$ square-matrix of elements $e_{jk} = (s+n_i)^{-1}(s+n_i+1)^{-1}(1+I_{\{j\}}(k))$ ($i=1$ and $2$, respectively), and $A$ is an $n_1\times n_2$ matrix with elements $a_{jk}=I_{(X_j,\infty)}(Y_k)$.
\end{Mtheorem}

\begin{Mcorollary}
\label{co:1}
For any $c\in[0,1]$, it holds that 
\begin{equation}
\label{eq:upmeandist}
\begin{array}{rcl}
&~& \overline{\mathcal{P}}\big[P(X\leq Y)> c|X^{n_1},Y^{n_2}\big]=P\left[g(w_{1\cdot},w_{2\cdot},X^{n_1}, Y^{n_2})>c\right],\\
\end{array}
\end{equation} 
with
$$
g(w_{1\cdot},w_{2\cdot},X^{n_1}, Y^{n_2})=w_{10}w_{20}+w_{10}\sum\limits_{j=1}^{n_2} w_{2j} +w_{20}\sum\limits_{j=1}^{n_1} w_{1j} +\sum\limits_{j=1}^{n_1}\sum\limits_{k=1}^{n_2} w_{1j}w_{2k} I_{(X_j,\infty)}(Y_k),
$$
where  $(w_{i0},w_{i1},\dots,w_{in_i})\sim Dir(s,\stackrel{n_i}{\overbrace{1,\dots,1}})$ for $i=1,2$.
Consider the augmented vectors $W=[w_{10},w_{11},\dots,w_{1n_1}]^T$, $V=[w_{20},w_{21},\dots,w_{2n_2}]^T$, and the matrix $A$ with elements $a_{jk} = I_{(X_{j-1},\infty)}(Y_{k-1})$ for all $j,k\neq{1}$ and $a_{jk} = 1$ if $j=1$ or $k=1$. The mean and variance of $g(w_{1\cdot},w_{2\cdot},X^{n_1}, Y^{n_2})$ can be computed using the same formulas as in (\ref{eq:lowmeanvar}), where, this time, $E[WW^T]$ and $E[VV^T]$ are $(n_i+1)\times (n_i+1)$ square-matrices ($i=1$ and $2$, respectively) of elements $e_{jk} = (s+n_i)^{-1}(s+n_i+1)^{-1}\tilde{e}_{jk}$ with $\tilde{e}_{jk}=(1+I_{\{j\}}(k))$ for all $j,k\neq{1}$ and $\tilde{e}_{jk} = s(1+s\delta_{jk})$ if $j=1$ or $k=1$.
\end{Mcorollary}
To perform the hypothesis test, we select $c=1/2$ and, according to the decision rule (\ref{eq:exploss_dec}) for some $K_0,~K_1$, we check if
$$
\underline{\mathcal{P}}\big[P(X\leq Y)> \tfrac{1}{2}|X^{n_1},Y^{n_2}\big]>1-\gamma, ~~\overline{\mathcal{P}}\big[P(X\leq Y)> \tfrac{1}{2}|X^{n_1},Y^{n_2}\big]>1-\gamma,
$$
where $\gamma=\tfrac{K_0}{K_0+K_1} \in (0,1)$ (e.g., $1-\gamma=0.95$).

\subsection{Choice of the prior strength $s$}
\label{sec:s}
The value of $s$ determines how quickly lower and upper posterior expectations converge at the increase
of the number of observations.
A way to select a  value of $s$ is by imposing that
the degree of robustness (indeterminacy) $\overline{\mathcal{E}}[P(X\leq Y)|X^{n_1},Y^{n_2}]-\underline{\mathcal{E}}[P(X\leq Y)|X^{n_1},Y^{n_2}]$
is reduced to a fraction of its prior value ($\overline{\mathcal{E}}[P(X\leq Y)]-\underline{\mathcal{E}}[P(X\leq Y)]=1$)
after one observation $(X_1,Y_1)$. % with $X_1<Y_1$.
Imposing a degree of imprecision close to $1$ after the first observation increases the probability of an indeterminate outcome of the test, whereas, a value close to $0$ makes the test less reliable (in fact the limiting value of $0$ corresponds to the BB which will be shown in  Section \ref{sec:sim} to be less reliable than the IDP). Then, the intermediate value of $1/2$ is a frequent choice in prior-ignorance modeling  \citep{pericchi_walley_91,walley1996}. Although this is a subjective way to choose the degree of conservativeness (indeterminacy), 
we will show in Section \ref{sec:sim} that it represents a reasonable trade-off between the reliability and indeterminacy of the decision.
From (\ref{eq:ranktest})  for $n_1=n_2=1$, it follows that
$$
\begin{array}{l}
\overline{\mathcal{E}}[P(X\leq Y)|X_1,Y_1]-\underline{\mathcal{E}}[P(X\leq Y)|X_1,Y_1]=\frac{s^2+2s}{(s+1)^2}. \\
\end{array}
$$
Thus,  by imposing that,
$$
\frac{s^2+2s}{(s+1)^2}=\frac{1}{2},
$$
we obtain $s= \sqrt{2}-1$. Observe that the lower and upper probabilities produced by a value of $s$ are always contained in the probability intervals produced by the larger value of $s$. Then, whenever we are undecided for $s_1$ we are also for $s_2 > s_1$. Nonetheless, for
large $n$ the distance between the upper and lower probabilities goes to 0, then also the indeterminateness  goes to zero.

\subsection{Managing ties}
\label{sec:ties}
To account for the presence of ties between samples from the two populations ($X_i=Y_j$), the common approach is to test the hypothesis $[P(X< Y) +\frac{1}{2} P(X=Y)]\leq 0.5$ against $[P(X< Y) +\frac{1}{2} P(X=Y)] > 0.5$. 
Since
$$
P(X< Y) +\frac{1}{2} P(X=Y)=E\left[I_{(X,\infty)}(Y)+\tfrac{1}{2}I_{\{X\}}(Y)\right]=E[H(Y-X)],
$$
where $H(\cdot)$ denotes the Heaviside step function, i.e., $H(z)=1$ for $z>0$,  $H(z)=0.5$ for $z=0$ and  $H(z)=0$ for $z<0$, in case of ties the $U$ statistic becomes
 \begin{equation}
\label{eq:Uties}
U=\sum\limits_{i=1}^{n_1} \sum\limits_{j=1}^{n_2} H(Y_j-X_i),
\end{equation}
and it represents the number of pairs $(X_i,Y_j)$ for which $X_i< Y_j$ plus half of the number of pairs $(X_i,Y_j)$ for which $X_i= Y_j$. 
The results presented in Section \ref{sec:MWW} are still valid if we substitute $I_{(X_j,\infty)}(Y_k)$ with $H(Y_k-X_j)$ in matrix $A$. 
%Also the convergence of the distribution of $P(X< Y) +\frac{1}{2} P(X=Y)$ with the distribution of the $U$ statistic under the null hypothesis can be extended to this case by computing the mean $\mu$ and variance $\sigma$ as in (\ref{eq:lowmeanvar}) with the new matrix $A$ and comparing it with the limiting mean and variance of $U$ for large $n$ in case of ties:
%$$
%\mu \approx 1/2, ~~~~ \sigma^2 \approx \frac{1}{6n}-\frac{\sum_{i=1}^k (t_i^3-t_i)}{12n^4}
%$$
%where, said $Z_1,\dots,Z_k$ the different values observed in $X$ and $Y$, $t_i$ is the number of observation in $X$ or $Y$ equal to $Z_i$.

\section{Asymptotic consistency}
\label{sec:Asymptotic}
From the expression of the lower and upper means in (\ref{eq:ranktest}), it can be verified that for $n_1,n_2\rightarrow\infty$:
$$
\underline{\mathcal{E}}[P(X\leq Y)|X^{n_1},Y^{n_2}],\overline{\mathcal{E}}[P(X\leq Y)|X^{n_1},Y^{n_2}] \sim {\mathcal{E}}[P(X\leq Y)|X^{n_1},Y^{n_2}]  \sim \frac{U}{n_1n_2}.
$$
Notice that in this section the symbol $\sim$ will be used to indicate asymptotic equivalence. 
The imprecision (degree of robustness) goes to zero for  $n_1,n_2\rightarrow \infty$ and the expectation ${\mathcal{E}}[P(X\leq Y)|X^{n_1},Y^{n_2}]$
is asymptotically equivalent to the Mann-Whitney statistic  \citep{Ferguson1973}.
The consistency of the IDP rank-sum test can be verified by considering the asymptotic behavior of the posterior lower and upper distributions of $P(X\leq Y)$ and compare it to the asymptotic distribution of the statistic $U/n_1n_2$. For ease of presentation, we limit ourselves to the case $n_1=n_2=n$.
In  \citet[Appendix A.5]{Lehmann} it is proved that $U/n_1n_2$ converges for $n_1,n_2\rightarrow\infty$ to a normal with mean $E[U_{ij}]=P(X\leq Y)$ and variance
\begin{equation} \label{eq:Uvar}
\frac{1}{n}Cov[U_{ij}, U_{i,k\neq j}]+\frac{1}{n}Cov[U_{ij},U_{k\neq i,j}],
\end{equation}
where $U_{rt}=I_{(X_r,\infty)}(Y_t)$ .
In the following theorem an equivalent result is proved for the lower distribution of $P(X \leq Y)$ in the IDP rank-sum test 
\begin{Mtheorem}
\label{th:2}
Assume that $n_1=n_2=n$, for $n\rightarrow \infty$ the IDP rank-sum test lower distribution converges to a normal with mean $E[U_{ij}]=P(X\leq Y)$ and variance given  by Equation (\ref{eq:Uvar}).
\end{Mtheorem}

The above proof can  be easily generalized to the upper distribution (the terms due to $w_{10}$ and $w_{20}$ vanish asymptotically) and to the case $n_1 \neq n_2$ (following the same procedure as in \citet[Th. 9]{Lehmann}). Theorem \ref{th:2} proves that the (upper and lower) distribution of the IDP rank-sum test is asymptotically equivalent to the distribution of the statistic $U/n_1n_2$ and, thus, the IDP rank-sum test  is consistent as a test for $P(X \leq Y)$.
Conversely, the MWW test is only consistent in the case $P(X \leq Y)=0.5$ and $F_X=F_Y$ or $P(X \leq Y) \neq 0.5$ and $F_X \neq F_Y$,
while it is not consistent for $P(X \leq Y) = 0.5$ and $F_X \neq F_Y$. 
For instance if $X \sim N(0,1)$ and $Y \sim N(0,\sigma^2)$ with $\sigma^2>1$,
two Normal distributions with different variance, then $P(X\leq Y)=0.5$ but the distributions are different.
In this case, if we apply MWW test with a significance level $\gamma=0.05$, MWW will return the alternative hypothesis
in approximatively $8.7\%$ of the cases (for a large $\sigma^2$), see \citet[Sec. 25.5]{dasgupta2008asymptotic}.
This means that MWW is not calibrated as a test for $P(X\leq Y)=0.5$ and it is not powerful as a test for  $F_X(x) \neq F_Y (x)$.
Conversely, because of Theorem \ref{th:2}, our IDP test with $\gamma=0.05$ will return the alternative hypothesis
(asymptotically) in  $5\%$ of the cases, which is correct since $P(X\leq Y)=0.5$.

\section{Numerical simulations}
\label{sec:sim}
Consider a Monte Carlo experiment in which $n_1,~n_2$ observations $X,Y$ are generated based  on
$$
X\sim N(0,1),~~Y\sim N(\Delta,1),
$$
with $\Delta$ ranging from $-1.5$ to $1.5$.
To facilitate the comparison of IDP tests with more traditional tests (which never issue indeterminate outcomes) we introduce a new test (called ``50/50 when indeterminate'') which returns the same response as the IDP when this is determinate,
and issues a random answer (with 50/50 chance) otherwise. We want to stress that the test ``50/50 when indeterminate'' has been introduced only for the sake of  comparison. We are not suggesting that when the IDP is indeterminate we should toss a coin to take the decision. On the contrary we claim that the indeterminacy of the IDP is an additional useful information that our approach gives to the analyst.  In these cases she/he knows that (i) her/his posterior decisions would depend on the choice of the prior $G_0$; (ii) deciding between the two hypotheses under test is a difficult problem as shown by the comparison with the Bayesian Bootstrap DP (BB-DP) rank-sum test ($s=0$) and MWW tests. Based on this additional information, the analyst can for example decide to collect additional measurements to  eliminate the indeterminacy (in fact we have seen that when the number of observations goes to infinity the indeterminacy goes to zero).
% We point the reader to \citet{zaffalon2012evaluating} for a broader discussion about indeterminacy and risk-aversion.  

\begin{figure}[h]
 \centerline{\includegraphics[width=12cm]{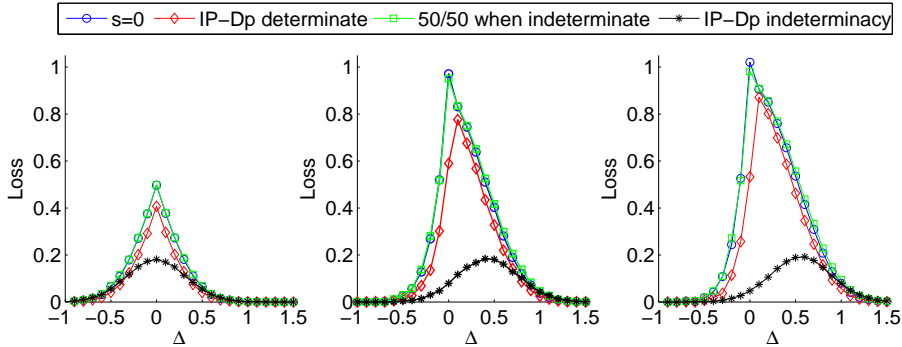}}
    \caption{Loss as a function of $\Delta$ for the case $K_0=K_1=1$ (left),  $K_0=1,K_1=9$ (center) and $K_0=1,K_1=19$ (right).}
\label{fig:1bsim}
\end{figure}

We start by comparing the performance of the BB-DP and IDP tests.
To evaluate the performance of the tests, we have used the loss function defined in (\ref{eq:loss}).
In particular, for each value of $\Delta$ we have performed  $20000$ Monte Carlo runs by generating in each run $n_1=n_2=20$
observations for $X,Y$. 
The average  loss for the cases (i) $K_1=K_2=1$ (i.e., $\gamma=0.5$) (ii) $K_1=1$ and $K_2=9$ (i.e., $\gamma=0.1$) and (iii) $K_1=1$ and $K_2=19$ (i.e., $\gamma=0.05$) is shown in Figure \ref{fig:1bsim} as a function of $\Delta$.
%\begin{figure}[htb]
 %\centerline{\includegraphics[height=4.5cm]{Figure_con_s=0_alpha100=50}}
    %\caption{Loss as a function of $\Delta$ for the case $K_0=K_1=1$.}
%\label{fig:1bsim}
%\end{figure}
%\begin{figure}[htb]
 %\centerline{\includegraphics[height=4.5cm]{Figure_con_s=0_alpha100=10}\includegraphics[height=4.5cm]{Figure_con_s=0_alpha100=5_fifty}}
    %\caption{Loss as a function of $\Delta$ for the case $K_0=1,K_1=9$ (left) $K_0=1,K_1=19$ (right).}
%\label{fig:1bsimb}
%\end{figure}
In particular, we report (i) the loss of the BB-DP test; (ii) the loss of the IDP test when it is determinate;
(iii) the indeterminacy of the IDP test, i.e., the number of times it returns an indeterminate response divided by the total number of Monte Carlo runs;
(iv) the loss of the ``50/50 when indeterminate'' test.

From  Figure \ref{fig:1bsim}, it is evident that the performance of the BB-DP and 50/50 tests practically coincide.
Furthermore, since we noticed from experimental evidence that in all cases in which BB-DP is determinate,  BB-DP returns the same response as IDP, 
%(i.e.,in the determinate cases the accuracy of BB-DP is also equal to the diamond curve in Figure \ref{fig:1bsim})
the difference between the two tests is only in the runs  where the IDP is indeterminate. In these runs, BB-DP is clearly guessing at random, since overall it has the same accuracy as the 50/50 test. Therefore,  the IDP is able to isolate several instances in which BB-DP is guessing at random, thus providing useful information to the analyst. 
Assume, for instance, that we are trying to compare the effects of two  medical treatments (``Y  is better than X'')
and that, given the available data, the IDP is indeterminate. In such situation the BB-DP test always issues a determinate response (I can tell if ``Y  is better than X''), but it turns out that its response is virtually random (like if we were tossing a coin). 
On the other side, the IDP acknowledges the impossibility of making a decision and thus, although BB-DP and the IDP (more precisely the ``50/50 when indeterminate'' test) have the same accuracy, the IDP provides more information. 
%The conclusions for the cases $K_1=1$ and $K_2=9$ (or $K_2=19$) are similar (see Figure \ref{fig:1bsim} center and right), in the sense that the performance of the BB-DP and 50/50 tests (almost) coincide.
Note that, in all the cases the maximum percentage of runs in which the IDP is indeterminate is about $18\%$;  this means that 
BB-DP is issuing a random answer in $18\%$ of the cases, which is a large percentage.
For large $|\Delta|$, i.e. when the hypothesis test is easy, there are not indeterminate instances and both the BB-DP and the IDP tests have zero loss.
It is interesting to note that,  for the cases $K_1=1$ and $K_2=9$ (or $K_2=19$) (Figure \ref{fig:1bsim} center and right) it is more risky (we may incur a greater loss) taking the action $a=1$ than $a=0$, and thus the indeterminacy curve is shifted to the $\Delta>0$ quadrant.

\begin{figure}[htb]
 \centerline{\includegraphics[height=4.5cm]{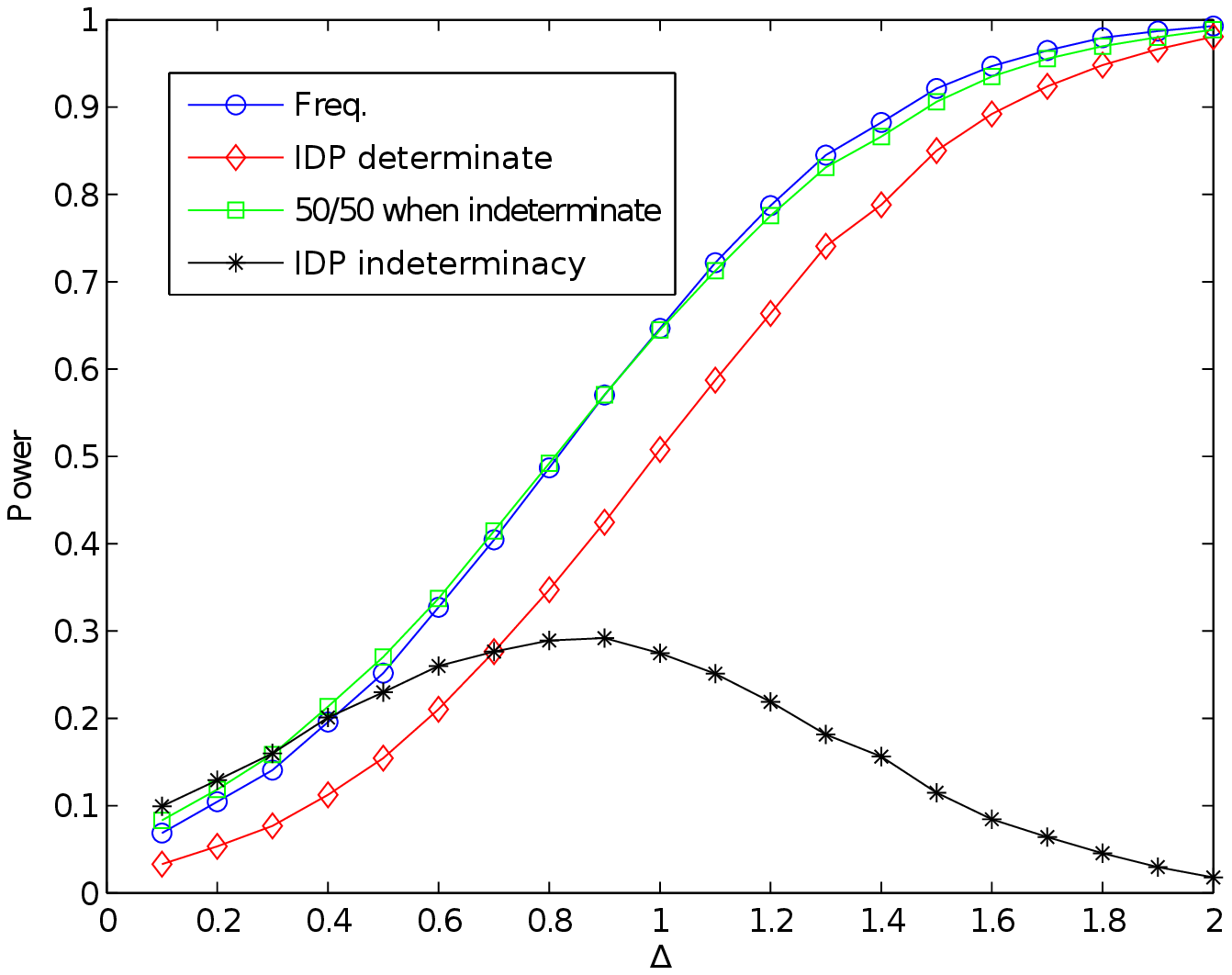}\includegraphics[height=4.5cm]{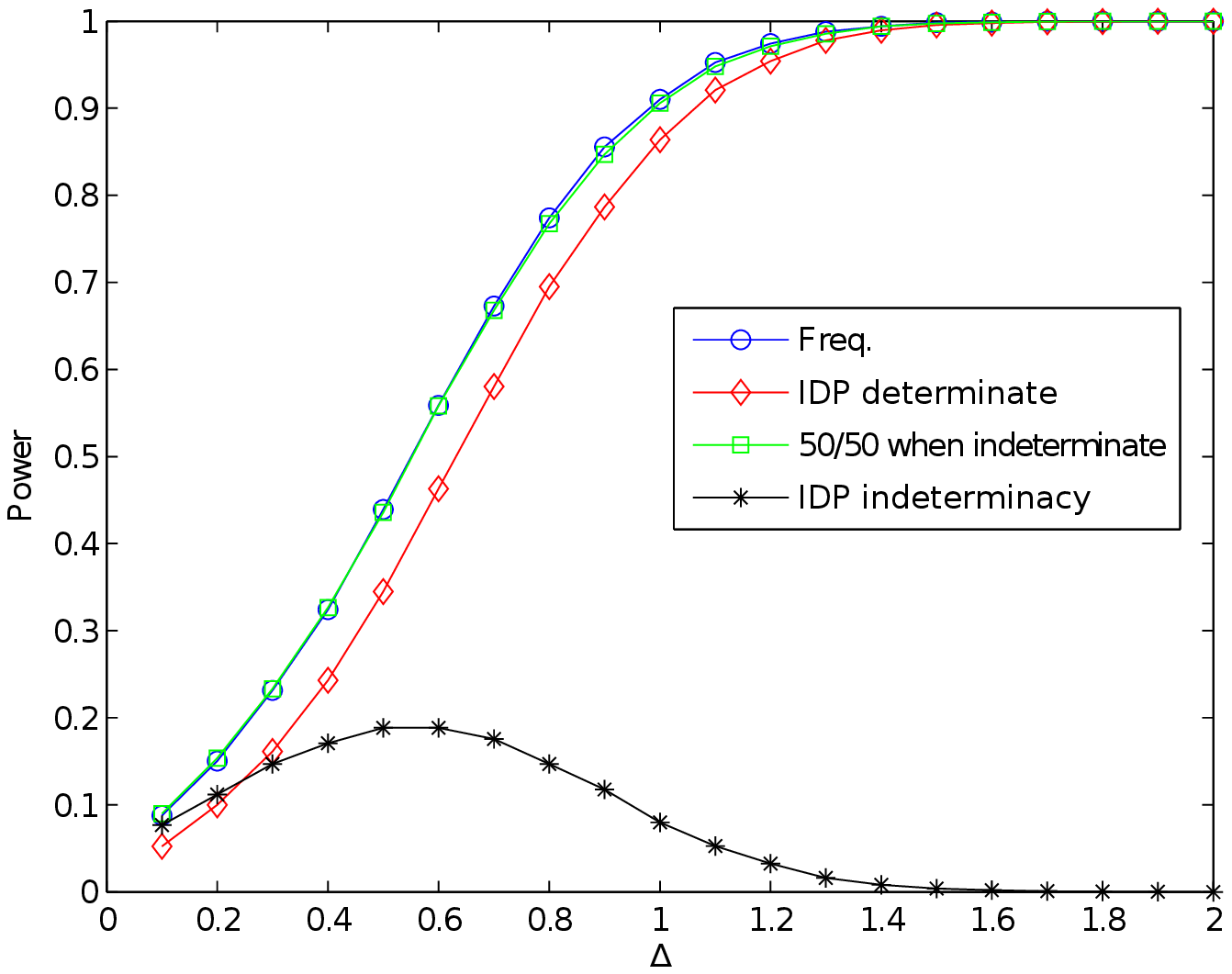}}
    \caption{Power as a function of the difference of the medians $\Delta$ for the case $n_1=n_2=10$ (left) and $n_1=n_2=20$ (right) with $\gamma=0.05$.}
\label{fig:1sim}
\end{figure}

\begin{figure}[htb]
  \centerline{\includegraphics[height=4.3cm]{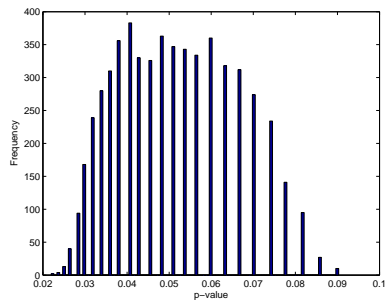}}
    \caption{Distribution of MWW p-values in the IDP indeterminate  cases  for $n_1=n_2=20$, $\gamma=0.05$ and $\Delta=0.5$.}
\label{fig:pval}
\end{figure}

We have also compared the IDP test and the one-sided MWW null hypothesis significance test (NHST) implemented  according to the conventional decision criterion,  $p<0.05$.  It is well known that the decision process in NHST is flawed. It is  based on asking
what is the probability of the data statistic if the null hypothesis were true. 
This means that NHST can only reject the null hypothesis ($\Delta \leq 0$), contrarily to a Bayesian analysis that can also accept this hypothesis.
Furthermore, in a Bayesian analysis we have a principled way to determine $\gamma$ (i.e., by means of a loss function) which is lost when putting decisions in the format $p < 0.05$ (or the more vague  $p < 0.1$). 
Because of these differences, it is difficult to compare the Bayesian with the NHST approach, where we do not have a clear interpretation of the significance level.
However, we believe a relatively fair comparison can be carried out by setting $\gamma$ equal to the significance level of the NHST test, so that the decision criteria adopted by the two test are as  similar as possible. 
Figure \ref{fig:1sim} shows the power for the case $\gamma=0.05$, $n_1=n_2=10$ and $n_1=n_2=20$.
In case $n_1=n_2=20$ (Figure \ref{fig:1sim}, right) it is evident that the performance of the MWW and 50/50 tests practically coincide.
Since it can be verified experimentally that when the IDP is determinate the two tests return the same results, this again suggests that when the IDP is indeterminate we have equal probability that $p < 0.05$ or  $p > 0.05$, as it is shown in Figure \ref{fig:pval}. The IDP test is able to isolate some instances in which also the MWW test is issuing a random answer.
Note that, for $\Delta=0.5$, the maximum percentage of runs in which the IDP test is indeterminate is large, about $18\%$;  this means that 
MWW is issuing a random answer in $18\%$ of the cases. 

% Looking at the curve of the IP-Dp indeterminacy, it can be observed that it reaches its maximum value
% at $\Delta\approx 0.6$ and then goes to zero for $\Delta > 1.5$.
% For $\Delta=2$ the accuracy of MWW and IP-Dp determinate is one, which means that when the hypothesis test
% is easy (large $\Delta$), there are not indeterminate instances and both tests have accuracy one.
The results for the case $n_1=n_2=10$ (Figure \ref{fig:1sim}, left) lead to similar conclusions. The performance of the MWW and 50/50 tests (almost) coincide.
The 50/50 test is slightly better for $\Delta\leq 0.9$ and slightly worse for $\Delta> 0.9$.
$\Delta=0.9$ is the value  which corresponds to the maximum indeterminacy of the IDP, i.e. $30\%$.
Thus, for $\Delta=0.9$, MWW is guessing at random in $30\%$ of the runs.

It is worth analyzing also the case $\Delta=0$. We know that in this case the frequentist test is calibrated, i.e., when $\gamma=0.05$
the percentage of correct answers is $95\%$ (although it can be noticeably larger for small values of $n_1,~n_2$ since the discreteness of the MWW statistic originates a gap between the chosen $\gamma$ and the actual significance of the MWW test).
Table \ref{tab:1sim} shows the accuracy for $\Delta=0$. The performance of the MWW and 50/50 tests are similar also in this case.
The difference is about $1\%$ (for $n_1=n_2=10$) and $0.5\%$ (for $n_1=n_2=20$). 
%Observe that, the accuracy $0.945$ of the 50/50 test is  obviously equal (on average) to  $0.911+0.068/2$, i.e.,
%the accuracy of IDP in the determinate cases plus half of the indeterminate cases.

\begin{figure}[htb]
  \centerline{\includegraphics[height=4.5cm]{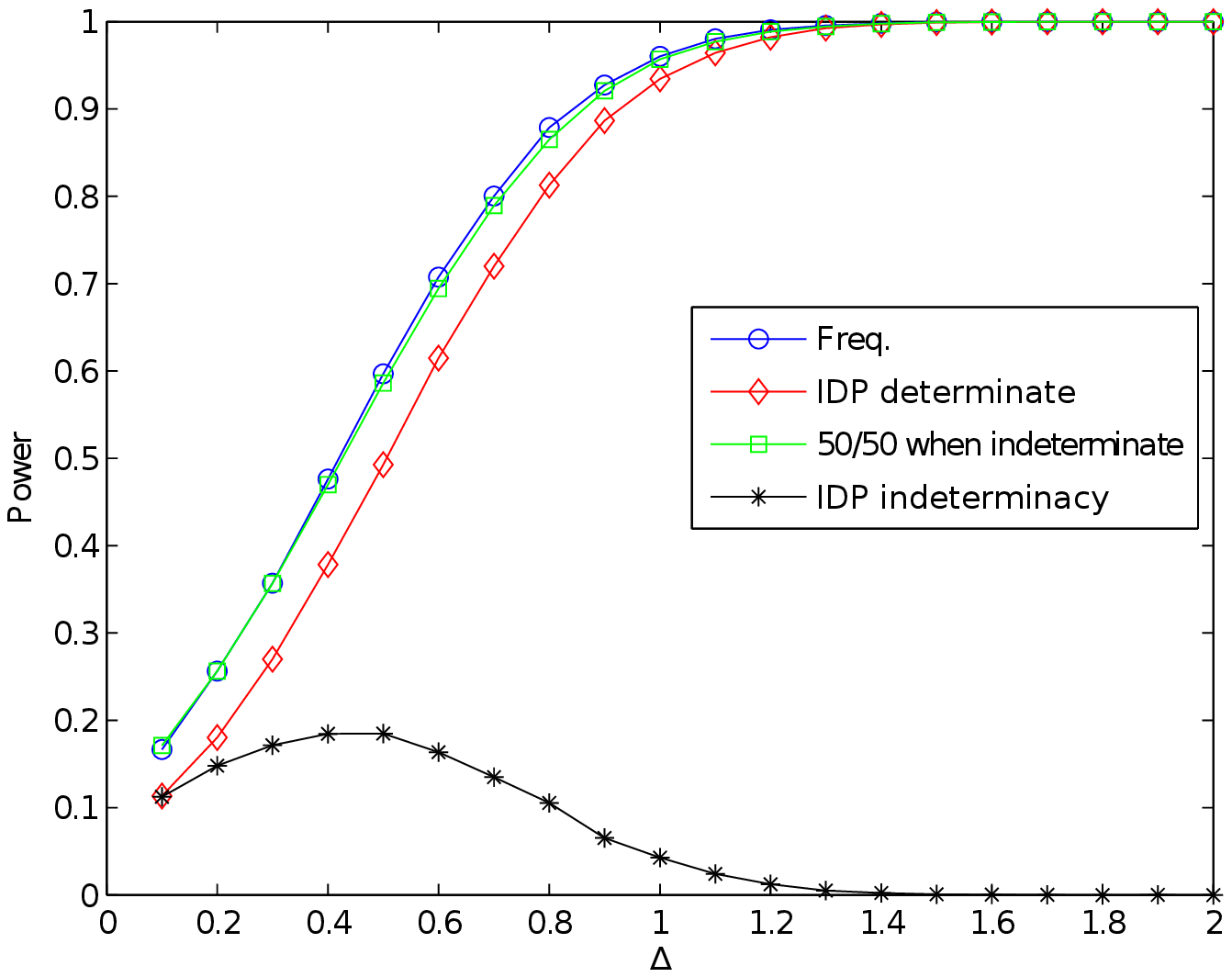}\includegraphics[height=4.5cm]{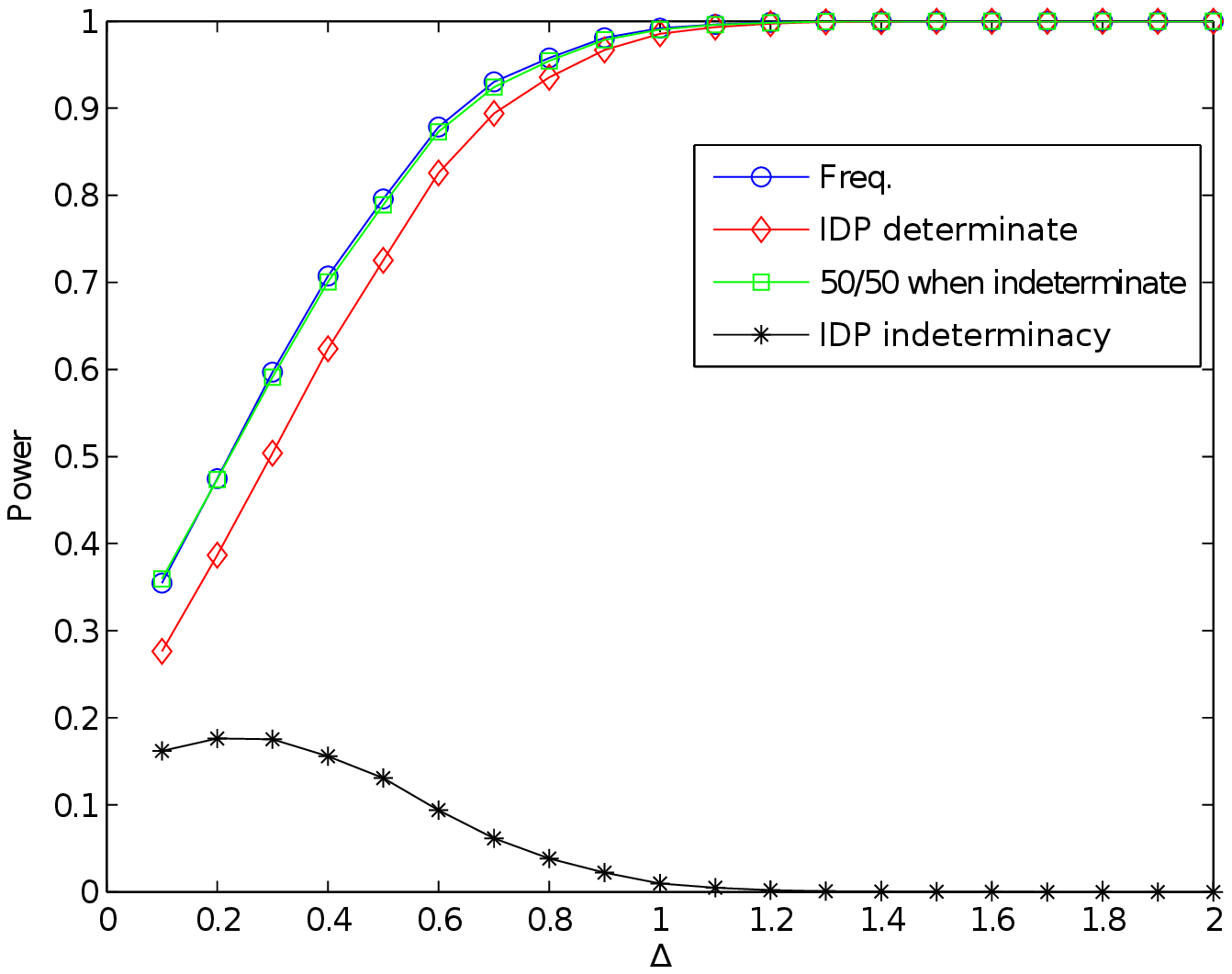}}
    \caption{Power as a function of the difference of the medians $\Delta$ for $n_1=n_2=20$, $\gamma=0.1$ (left) and $\gamma=0.25$ (right).}
\label{fig:2sim}
\end{figure}

\begin{table}
\centering
{\small
 \begin{tabular}{|c|c|c|c}
\hline
& Accuracy $n_1=n_2=10$ & Accuracy $n_1=n_2=20$\\
\hline
\hline
MWW & 0.955 & 0.952\\
50/50 test & 0.945 & 0.947\\
\hline
IDP when determinate & 0.911 & 0.924\\
Indeterminacy & 0.068 & 0.045\\
\hline
 \end{tabular}  }
\caption{Accuracy for  $\Delta=0$ and $\gamma=0.05$.}
\label{tab:1sim}
\end{table}
Also in this case, when the IDP is determinate, it returns the same responses as MWW.  This result holds independently of the choice of $\gamma$, as shown by Figure \ref{fig:2sim} and Table \ref{tab:2sim} where we have repeated the above experiment for $n_1=n_2=20$ with, this time, $\gamma=0.1$ and $\gamma=0.25$.

\begin{figure}[htb]
  \centerline{\includegraphics[height=4.5cm]{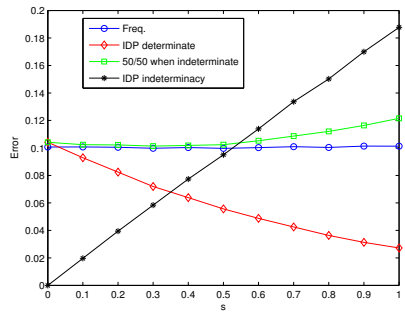}}
    \caption{Error as a function of $s$ for $n_1=n_2=20$, $\gamma=0.1$ and $\Delta=0$.}
\label{fig:2simb}
\end{figure}

Finally, Figure \ref{fig:2simb} shows the error (one minus the accuracy) of the IDP test as a function of $s$, when $\gamma=0.1$, $n_1=n_2=20$ and $\Delta=0$.
Clearly, the error of the MWW test is constantly equal to $\gamma=0.1$ (we are under the null hypothesis of MWW).
The error of the IDP test when determinate decreases with $s$, because of the increase of the indeterminacy.
%Note that for $s=0$, the error of the IDP when determinate and the error of the 50/50 test coincide, since there are not indeterminate instances in this case.
The error of the 50/50 test has a convex trend, clearly decreasing for $ s< 0.2$ and increasing for $s>0.5$. 
%It is higher for $s=0$ (Bayesian Boostrap)  than for  and then it increases again for $s>0.5$.
This (together with the other results of this section) may be seen as an empirical confirmation that the choice of $s=\sqrt{2}-1$ is appropriate, since it guarantees a good trade-off between robustness and indeterminacy.

Finally observe that all the above differences/similarities between the three tests appear also 
in the case we consider location-shift models with distributions different from Gaussians (e.g., 
Student-t distribution). These results have been omitted for shortness.

\begin{table}
\centering
{\small
 \begin{tabular}{|c|c|c|c}
\hline
& Accuracy $\gamma=0.1$ & Accuracy $\gamma=0.25$ \\
\hline
\hline
MWW & 0.8995 & 0.7552\\
50/50 test & 0.8993 & 0.7482\\
\hline
IDP when determinate & 0.8568 & 0.6777\\
IDP indeterminacy & 0.081 & 0.142\\
\hline
 \end{tabular} }
\caption{Accuracy in case $\Delta=0$ for $n_1=n_2=20$}
\label{tab:2sim}
\end{table}

\section{Conclusions}
In this paper we have proposed a model of prior ignorance for nonparametric inference based on the Dirichlet process (DP), by extending the approach proposed by \citet{pericchi_walley_91} and based on the use of sets of prior distributions. 
We developed a prior near-ignorance DP model (IDP) for inference about a variable $X$ by fixing the prior strength of the DP and  letting the normalized probability measure vary in the set of all distributions. We have proved that the IDP is consistent and a-priori vacuous for all predictive inferences that can be defined as the expectation of a real-valued bounded function of $X$. 
The proposed IDP model has two main merits. First, it removes the need for specifying the infinite-dimensional parameter of the DP (only an upper bound for the strength $s$ of the DP must be assumed a-priori), thus making the elicitation of the prior very easy. Second, it allows computing the posterior inferences for which no closed form expression exists, by a simple Monte Carlo sampling 
from the Dirichlet distribution, thus avoiding more demanding sampling approaches typically used for the DP (e.g., stick breaking). 
Based on this new prior near-ignorance model, we have proposed a general, simple and conservative approach to Bayesian nonparametric tests, and in particular we have developed a robust Bayesian alternative to the Mann-Whitney-Wilcoxon test: the IDP rank-sum test.
We have shown that our test is asymptotically consistent, while this is not always the case for the  Mann-Whitney-Wilcoxon test.
Finally, by means of numerical simulations, we have compared the IDP rank-sum test to the Mann-Whitney-Wilcoxon test and the Bayesian test obtained from the DP when the prior strength goes to zero. Results have shown that the IDP test is more robust, in the sense that it is able to isolate instances in which these tests are almost guessing at random.
% test outperforms both the frequentist Mann-Whitney-Wilcoxon test and  the Bayesian test obtained as a limit of DP when the prior strength goes to zero.
Given these interesting  results, as future work we plan to use this approach to implement Bayesian versions of the most used frequentist nonparametric tests. % such as, for instance, sign test, Wilcoxon signed-rank test, test for censored data, test for bivariate dependence etc.
In the long run, our aim is to build a statistical package for Bayesian nonparametric tests.

\section{Appendix}
{\it Proof of Theorem \ref{th:PI}:}  From (\ref{eq:expf}) assuming that $P\sim Dp(s,\alpha^*)$ one has that ${\mathcal{E}}[E(f)]=\int f d\alpha^*$.
Define $x_l= \arg\inf_{x \in \mathbb{X}} f(x)$ and $x_u= \arg\sup_{x \in \mathbb{X}} f(x)$, then (\ref{eq:vac1}) 
% it is enough to consider the DPs $Dp(s,\delta_{x_l})$ and $Dp(s,\delta_{x_u})$ belonging to the class (\ref{eq:IDP}), where
% $x_l= \arg\inf_{x \in \mathbb{X}} f(x)$ and $x_u= \arg\sup_{x \in \mathbb{X}} f(x)$:
follows by:
{\small
$$
\underline{\mathcal{E}}[E(f)]=\inf_{\alpha^* \in \mathbb{P}} \int f d\alpha^*= \int f d\delta_{x_l}=f(x_l),~~\overline{\mathcal{E}}[E(f)]=\sup_{\alpha^* \in \mathbb{P}} \int f d\alpha^*= \int f d\delta_{x_u}=f(x_u),
$$
}
which are the infimum and supremum of $f$ by definition. The lower and upper bounds are thus obtained by the following degenerate DPs $Dp(s,\delta_{x_l})$  and $Dp(s,\delta_{x_u})$, which belong to the class (\ref{eq:IDP}).
% In case $x_l$ (similar for $x_u$) does not belong to $\mathbb{X}$, the lower bound is not obtained by a DP with atomic base measure, but 
% $Dp(s,\delta_{x_l})$
In case $x_l$ is equal to $\infty$ (or $-\infty$), with $f(x_l)$ we mean $\lim_{x_l \rightarrow  \infty} f(x_l)$,  similar for the upper.\\
\noindent
{\it Proof of Theorem \ref{th:DPpost}:}  By exploiting the fact that  ${\mathcal{E}}[E(f)|X_1,\dots,X_n]=\int f ~d(\tfrac{s}{s+n}\alpha^*+\tfrac{n}{s+n}\frac{1}{n}\sum_{i=1}^n \delta_{X_i})$, the proof is similar to that of Theorem \ref{th:PI} (the lower and upper bounds are again obtained by degenerate DPs $Dp(s,\delta_{x_l})$ and $Dp(s,\delta_{x_u})$).\\
\noindent
{\it Proof of  Lemma \ref{ml:1}:} It follows from the properties (a) and (c) of the DP in Section \ref{sec:DP}.\\
% definition of Dirichlet process and the discreteness of the support of $G_{n_i}$, which imply the discreteness of the support of $F_{n_i}$ and the Dirichlet distribution for the vector of probabilities $(P(\{Z_{0}\}),P(\{Z_{1}\}),\dots, P(\{Z_{n_1}\}))$.\\
\noindent
{\it Proof of Theorem \ref{th:1}:} 
Based on the stick-breaking construction, a sample $F_0$ from the generic DP $Dp(s,G_0)$ can be written as
$F_0(x)=\sum\limits_{k=1}^{\infty} \pi_k \delta_{\tilde{X}_k}$ 
where $\pi_k=\beta_k\prod\limits_{i-1}^{k-1}(1-\beta_i)$, $ \beta_k \sim Beta(1,s)$, and $ X_k \sim G_0$. 
Then, using  (\ref{eq:mixing}), we have that
%Also, given the independent random variables $(w_0,w_1,\dots,w_n)\sim Dir(s,\stackrel{n_i}{\overbrace{1,\dots,1}})$ and distributions $P_0\sim Dp(s,G_0)$,$P_i\sim Dp(1,\delta_{X_i}$, the mixture
%$$
%\sum\limits_{i=0}^n w_i P_i
%$$
%has a DP distribution with base measure $s G_0 + \sum_{i=1}^n \delta_{X_i}$. Then, a sample from the posterior DP given the observations $X_1,\dots,X_n$ can be written as
\begin{equation} \label{eq:stick}
F_{n}(x)= \sum\limits_{i=1}^{n} w_i \delta_{X_i} + w_{0}\sum\limits_{k=1}^{\infty} \pi_{k} \delta_{\tilde{X}_k},
\end{equation}
where $(w_0,w_1,\dots,w_n)\sim Dir(s,\stackrel{n_i}{\overbrace{1,\dots,1}})$. 
Consider the two samples $F_X(x)$ and $F_Y(y)$ from the posterior distributions of $X$ and $Y$ given the generic DP priors $Dp(s,G_{10})$ and $Dp(s,G_{20})$.
%they can be written as
%\begin{equation} \label{eq:stick}
%\begin{array} {ll}
%F_{n_1}(x)=& \sum\limits_{i=1}^{n_1} w_{1i} \delta_{X_i} + w_{10}\sum\limits_{i=n_1+1}^{\infty} \pi_{1i} \delta_{\tilde{X}_i}\\
%F_{n_2}(y)=& \sum\limits_{j=1}^{n_2} w_{2j} \delta_{Y_j} + w_{20}\sum\limits_{j=n_2+1}^{\infty} \pi_{2j} \delta_{\tilde{X}_j}\\
%\end{array}
%\end{equation}
The probability of $P(X\leq Y)>c$ is $\mathcal{P}\big[P(X\leq Y)>c \big] = \mathcal{P}[\int F_{n_1}(y) d F_{n_2}(y)>c ]$.
Then, the posterior lower probability of $P(X\leq Y)> c$ is obtained by minimizing $\int F_{n_1}(y) d F_{n_2}(y)$, which, by  (\ref{eq:stick}), is equal to
\begin{equation} \label{eq:postP}
{\small
\begin{array}{ll}
~&\int \left(\sum\limits_{i=1}^{n_1} w_{1i} I_{(X_i,\infty)}(y)+w_{10} \sum\limits_{k=1}^{\infty} \pi_{1k} I_{(\tilde{X}_k,\infty)}(y)\right)
             \left(\sum\limits_{j=1}^{n_2} w_{2j} \delta_{Y_j}(y) + w_{20} \sum\limits_{l=1}^{\infty} \pi_{2l} \delta_{\tilde{Y}_l}(y)\right)dy\\
=&\sum\limits_{i=1}^{n_1}\sum\limits_{j=1}^{n_2} w_{1i} w_{2j} I_{(X_i,\infty)}(Y_j)
                   +w_{20} \sum\limits_{i=1}^{n_1}\sum\limits_{l=1}^{\infty} w_{1i} \pi_{2l} I_{(X_i,\infty)}(\tilde{Y}_l)\\
                  +&w_{10} \sum\limits_{k=1}^{\infty} \sum\limits_{j=1}^{n_2} \pi_{1k} w_{2j} I_{(\tilde{X}_k,\infty)}(Y_j)
									+w_{10}w_{20} \sum\limits_{k=1}^{\infty}  \sum\limits_{l=1}^{\infty} \pi_{1k} \pi_{2l} I_{(\tilde{X}_k,\infty)}(\tilde{Y}_{l})
\end{array}
}
\end{equation}
The minimum of $\int F_{n_1}(y) d F_{n_2}(y)$ is always found in correspondence of prior DPs such that the posterior probability of sampling $\tilde{X}_k<Y_j,\tilde{Y}_{l}$ or $\tilde{Y}_{l}>X_i$ is zero, so that only the term $\sum\limits_{i=1}^{n_1}\sum\limits_{j=1}^{n_2} w_{1i} w_{2j} I_{(X_i,\infty)}(Y_j)$ remains in  (\ref{eq:postP}). Priors of such kind are, for example, the extreme DP priors that give the posterior lower mean in (\ref{eq:ranktest}) and the posterior Dirichlet process $F_{n_i}$ with base probability $G_{n_i}$ given by (\ref{eq:basePost}).
%First, observe that the posterior lower probability of $P(X\leq Y)> a$ is obtained in correspondence of the extreme DP prior that gives the posterior lower mean in (\ref{eq:ranktest}) and the posterior Dirichlet process $F_{n_i}$ with base probability $G_{n_i}$ given by (\ref{eq:basePost}). Then the lower probability of $P(X\leq Y)>a$ is:
%$$
%\mathcal{P}\big[P(X\leq Y)>a|X^{n_1},Y^{n_2}\big] = \mathcal{P}\left[\displaystyle{\int F_{n_1}(y) d F_{n_2}(y)>a}\right],
%$$
%which, by Lemma \ref{ml:1}, is equal to
%$$
%\begin{array}{ll}
%~&\mathcal{P}\left[\displaystyle{\int \left(w_{10}I_{(X_0,\infty)}(y)+\sum_{j=1}^{n_1} w_{1j} I_{(X_j,\infty)}(y)\right)\left(w_{20}\delta_{Y_0}(y)+\sum_{j=1}^{n_2} w_{2j} \delta_{Y_j}(y)\right)dy}>a\right]\\
%=&\mathcal{P}\left[\sum\limits_{j=1}^{n_1}\sum\limits_{k=1}^{n_2} w_{1j}w_{2k} I_{(X_j,\infty)}(Y_k)>a\right].
%\end{array}
%$$
%
From the property of the Dirichlet distribution, we know that $E[w_{ij}]=1/(s+n_i)$ and, thus, we can rewrite the lower expectation given in the first equation of (\ref{eq:ranktest}) as
$$
\mu= \sum\limits_{j=1}^{n_1}\sum\limits_{k=1}^{n_2} \frac{1}{s+n_1} \frac{1}{s+n_2} I_{(X_j,\infty)}(Y_k)=E_W[W]^TAE_V[V],
$$
For the variance, we have that $\sigma^2=E[ (\sum\limits_{j=1}^{n_1}\sum\limits_{k=1}^{n_2}w_{1j}w_{2k} I_{(X_j,\infty)}(Y_k) )^2]-\mu^2$. Thus, by exploiting the equality
$$
\left(\sum\limits_{j=1}^{n_1}\sum\limits_{k=1}^{n_2}w_{1j}w_{2k}  I_{(X_j,\infty)}(Y_k) \right)^2=W^TAVW^TAV=V^TA^TWW^TAV,
$$
the linearity of expectation and the independence of $W,V$, one obtains
$$
E[V^TA^TWW^TAV]=E_V[V^TA^TE_W[WW^T]AV]=E_V[V^TA^TE_W[WW^T]AV].
$$
Since the result of this product is a scalar, it is equal to its trace and thus we can use the cyclic property  
$trace [E_V[V^TA^TE_W[WW^T]AV]] = trace[A^TE_W[WW^T]AE_V[VV^T]]$, and finally obtain
$\sigma^2= trace[A^TE_W[WW^T]AE_V[VV^T]]-\mu^2$. The proof is easily completed by deriving $E_W[WW^T]$ and $E_V[VV^T]$ from the fact that $w_{ij},w_{kl}$ are independent and
$E_W[w^2_{ij}]=\tfrac{2}{(s+n_i)(s+n_i+1)}$, $E_W[w_{ij}w_{il}]=\tfrac{1}{(s+n_i)(s+n_i+1)}$.\\
% In matrix form, the variance can be rewritten as
% $$
% \begin{array}{rcl}
% \sigma^2&=& trace\left( A^TE_W[WW^T]AE_V[VV^T]- A^TE[W]E[W]^TAE[V]E[V]^T \right),\\
% &=& tr\left( A^TE_W[WW^T]AE_V[VV^T]- A^TE[W]E[W]^TAE[V]E[V]^T \right),\\
% \end{array}
% $$
\noindent
{\it Proof  of Corollary \ref{co:1}:} 
First, observe that the posterior upper probability of $P(X\leq Y)> c$ is obtained in correspondence of the extreme DP prior that gives the posterior upper mean in (\ref{eq:ranktest}) and has base probability $dG_{n_i}$ given by (\ref{eq:basePost}). The probability of $X\leq Y$ for a given realization $F_{n_1}$ of $Dp(s,G_{n_1})$, and $F_{n_2}$ of $Dp(s,G_{n_2})$ is: 
$$
\begin{array}{ll}
~&\mathcal{P}\big[P(X\leq Y)>c|X^{n_1},Y^{n_2}\big] = \mathcal{P}\left[\displaystyle{\int F_{n_1}(y) d F_{n_2}(y)}>c\right] \vspace{2mm}\\
=&\mathcal{P}\left[\displaystyle{\int \left(w_{10}I_{(X_0,\infty)}(y)+\sum_{j=1}^{n_1} w_{1j} I_{(X_j,\infty)}(y)\right)\left(w_{20}\delta_{Y_0}(y)+\sum_{j=1}^{n_2} w_{2j} \delta_{Y_j}(y)\right)dy}>c\right]\\
=&\mathcal{P} \left[ w_{10}w_{20}+w_{10}\sum\limits_{j=1}^{n_2} w_{2j} +w_{20}\sum\limits_{j=1}^{n_1} w_{1j} +\sum\limits_{j=1}^{n_1}\sum\limits_{k=1}^{n_2} w_{1j}w_{2k} I_{(X_j,\infty)}(Y_k)>c\right].
\end{array}
$$
The computations are similar to those in Theorem \ref{th:1}, but in this case we must also consider 
the expectations $E_W[w^2_{i0}]=s(s+1)/(s+n_i)(s+n_i+1)$, $E_W[w_{i0}w_{ij}]=s/(s+n_i)(s+n_i+1)$ for $j>0$.\\
\noindent
{\it Proof of Theorem \ref{th:2}:} 
Our goal is to prove the convergence to a normal distribution of the Bayesian bootstrapped two-sample statistic 
$U_{DP} = \sum\limits_{i,j} w_{1i}w_{2j}I_{[X_i,\infty)}(Y_j)$, which implies the asymptotic normality of the DP rank sum test lower distribution, since the contribution of the prior $G_0$ vanishes asymptotically.
The asymptotic normality of $U_{DP}$ can be proved by means of Lemma 6.1.3. of \citet{Lehmann1998}, which states that given a sequence of random variables $T_n$, the distributions of which tend to a limit distribution $L$, the distribution of another sequence $T_n^*$ satisfying $E[(T_n^*-T_n)^2]\rightarrow 0$ also tends to $L$.  
Said $h(x,y) = I_{[x,\infty)}(y)$, $h_1(x) = E_Y[h(x,Y)]$ and  $h_2(y) = E_X[h(X,y)]$, the theorem will be proved by applying the lemma to 
$$
T_n = \sqrt{n}\left[ \frac{1}{n}(\sum_{i=1}^n h_1(X_i) - \theta) +\frac{1}{n}(\sum_{j=1}^n h_2(Y_j)-\theta)\right]
$$
and $T_n^* = \sqrt{n}(U_{DP}-\theta)$ where $\theta = E[U_{ij}] =  E[h_1(X)] = E[h_2(Y)] $.
$T_n$ is a sum of independent terms and thus, from the central limit theorem, it converges to a Gaussian distribution with mean $0$ and variance $\sigma^2=\sigma^2_1+\sigma^2_2$, where $\sigma^2_1=Var[h_1(X)]$ and $ \sigma^2_2=Var[h_2(Y)])$. Note that 
$$
\begin{array} {l}
\sigma_1^2 = Cov[h(X,Y),h(X,Y')] = Cov[U_{ij}, U_{i,k\neq j}],\\
\sigma_2^2 = Cov[h(X,Y),h(X',Y)] = Cov[U_{ij}, U_{k\neq i,j}].\\
\end{array}
$$
From Theorem \ref{th:1}, the mean of the lower distribution of $U_{DP}$ is $ \mu_l=E_W[W]^TAE_V[V]=\tfrac{U}{(s+n)^2}$, 
and thus, for large $n$, it is asymptotic to $U/n^2$ which converges, in turn, to $E[U_{ij}] = \theta$.%, see \citet[Appendix A.5]{Lehmann}
Then, also $E[T_n^*] = 0$ so that 
$$
E[(T_n^*-T_n)^2] = Var[T_n^*]+Var[T_n]-2 Cov[T_n^*,T_n].
$$
The proof will be completed by showing that $Var[T_n^*] \rightarrow \sigma^2$ and $Cov[T_n^*,T_n] \rightarrow \sigma^2$.
% $$
% \begin{array}{cc}
% Var[T_n^*] \rightarrow \sigma^2, & Cov[T_n^*,T_n] \rightarrow \sigma^2
% \end{array}
% $$
For the variance of $U_{DP}$ (\ref{eq:lowmeanvar}), first note that we can rewrite
$E_W[WW^T]=E_V[VV^T]=(D+J_n)\tfrac{1}{(s+n)(s+n+1)}$
where $D$ is the diagonal matrix of ones (identity matrix) and $J_n$ is the $n \times n$ matrix of ones.
Thus, we have that
$$
A^TE_W[WW^T]AE_V[VV^T]=A^T(D+J_n)A(D+J_n)\frac{1}{(s+n)^2(s+n+1)^2},
$$
and, for large $n$,
$$
\begin{array}{l}
\frac{trace(A^T(D+J_n)A(D+J_n))}{(s+n)^2(s+n+1)^2}\rightarrow \frac{trace(A^TA) +trace(A^TAJ_n)+trace(A^TJ_nA)+trace(A^TJ_nAJ_n)}{n^2(n+1)^2}.
\end{array}$$
The above sum has four terms at the numerator:
{\small
$$
\begin{array}{lcl}
trace(A^TA)&=&\sum_{i,j} a_{ij}^2=\sum_{i,j} I_{(X_i,\infty)}(Y_j),\\
trace(A^TAJ_n)&=&\sum_{i,j} a_{ij} \sum_{k} a_{ik}=\sum_{i,j} I_{(X_i,\infty)}(Y_j)+\sum_{i,j\neq k} I_{(X_i,\infty)}(Y_j)I_{(X_i,\infty)}(Y_k),\\
trace(A^TJ_nA)&=&\sum_{i,j} a_{ij} \sum_{k} a_{kj}=\sum_{i,j} I_{(X_i,\infty)}(Y_j)+\sum_{i\neq k,j} I_{(X_i,\infty)}(Y_j)I_{(X_k,\infty)}(Y_j),\end{array}
$$
}
and $
trace(A^TJ_nAJ_n)=trace(A^T\mathbbm{1}\mathbbm{1}^TA\mathbbm{1}\mathbbm{1}^T)=trace(\mathbbm{1}^TA\mathbbm{1}\mathbbm{1}^TA^T\mathbbm{1})=U^2,
$
where $\mathbbm{1}$ is the unit vector.
% Following \citet[Example 5, Appendix A.1]{Lehmann}, let us introduce 
% $$
% p_1=P(X\leq Y), ~~p_2=P(X\leq Y,X\leq Y'),~~p_3=P(X\leq Y,X'\leq Y),
% $$ 
Then we have that
$$
\begin{array}{l}
\sigma_l^2=
\frac{3\sum_{i,j} I^2_{(X_i,\infty)}(Y_j)-3n^2\mu_l^2 }{n^2(n+1)^2}
+\frac{\sum_{i,j\neq k} I_{(X_i,\infty)}(Y_j)I_{(X_i,\infty)}(Y_k)-n^2(n-1)\mu_l^2}{n^2(n+1)^2}\\
+\frac{\sum_{i\neq k,j} I_{(X_i,\infty)}(Y_j)I_{(X_k,\infty)}(Y_j)-n^2(n-1)\mu_l^2}{n^2(n+1)^2}
+\frac{3n^2+2n^2(n-1)+n^4}{n^2(n+1)^2}\mu_l^2-\mu_l^2.
\end{array}
$$
Note that $(\tfrac{3n^2+2n^2(n-1)+n^4}{n^2(n+1)^2}-1)\mu_l^2=0$ and, since the first term in $\sigma_l^2$ goes to zero as $1/n^2$, for large $n$,
$$
\begin{array}{l}
\sigma_l^2 \rightarrow 
% \frac{3\sum_{i,j} I^2_{(X_i,\infty)}(Y_j)-3n^2\mu_l^2 }{n^2(n+1)^2}\\
\frac{\sum_{i,j\neq k} I_{(X_i,\infty)}(Y_j)I_{(X_i,\infty)}(Y_k)-n^2(n-1)\mu_l^2}{n^2(n+1)^2}
+\frac{\sum_{i\neq k,j} I_{(X_i,\infty)}(Y_j)I_{(X_k,\infty)}(Y_j)-n^2(n-1)\mu_l^2}{n^2(n+1)^2}.
% +\frac{3n^2+2n^2(n-1)+n^4}{n^2(n+1)^2}\mu_l^2-\mu_l^2.
\end{array}
$$
For large $n$, it can be shown  \citet[Th. 9]{Lehmann} that the right-hand side of the above equations tends to
$\frac{1}{n}Cov[U_{ij}, U_{i,k\neq j}] +\frac{1}{n}Cov[U_{ij},U_{k\neq i,j}] = \frac{1}{n}\sigma^2$, and thus $Var[T_n]=Var[\sqrt{n}U_{DP}] \rightarrow \sigma^2$.
For the covariance we have
{\small
$$
\begin{array}{l}
Cov[T_n,T^*_n] = \left(E[U_{DP}\sum\limits_{i=1}^n h_1(X_i)]+E[U_{DP} \sum\limits_{j=1}^n h_2(Y_j)] -2\theta \right)\\
~= \left(E[\sum\limits_{i,j} w_{1i}w_{2j} E_{Y_j}[h(X_i,Y_j)]\sum\limits_{i=1}^n h_1(X_i)]+E[\sum\limits_{i,j} w_{1i}w_{2j} E_{X_i}[h(X_i,Y_j)] \sum\limits_{j=1}^n h_2(Y_j)] -2\theta \right)\\
~= \left(E_{X}[\sum\limits_{i=1}^n E[w_{1i}] h_1(X_i) \sum\limits_{i=1}^n h_1(X_i)]+E_{Y}[\sum\limits_{j=1}^n E[w_{2j}] h_2(Y_j) \sum\limits_{j=1}^n h_2(Y_j)] -2\theta \right)\\
~= \frac{1}{n} \left((\sum\limits_{i=1}^n E_{X}[h_1(X_i))^2]+(\sum\limits_{j=1}^nE_{Y}[ h_2(Y_j))^2] -2\theta \right)= Var[h_1(X)] + Var[h_2(Y)] = \sigma^2 .
\end{array}
$$}

%%%%%%%%%%%%%%%%%%%%%%%%%%%%%%%%%%%%%%%%%%%%%%%%%%%%%%%%%%%%%%%%%%%%%%%%%%%%%%%%%%%%%%%%%%%%%%%%%%%%%%%%%%%%%%%%%%%%%%%%%%%%

\section*{Acknowledgements} This work  has been partially supported by the Swiss NSF grant
n.~200020-137680/1.

%%%%%%%%%%%%%%%%%%%%%%%%%%%%%%%%%%%%%%%%%%%%%%%%%%%%%%%%%%%%%%%%%%%%%%%%%%%%%%%%%%%%%%%%%%%%%%%%%%%%%%%%%%%%%%%%%%%%%%%%%%%%

% \noindent{\large\bf References}
% \begin{description}
% \item
% Andrews, D. W. K. (1984). Non-strong mixing autoregressive processes.
%   {\it J. Appl. Probab.} {\bf 21}, 930-934.
% \item
% Avram, F. and Taqqu, M. S. (1987)  Noncentral limit theorems and Appell
% polynomials. {\it Ann. Probab.} {\bf 15}, 767-775.
% \item
% Bradley, R. C. (1986). Basic properties of strong mixing conditions. In {\it
% Dependence
% in Probability and Statistics} (Edited by E. Eberlein and M. S. Taqqu),
% 162-192. Birkh\"auser, Boston.
% \item
% Fox, R. and Taqqu, M. (1987). Central limit theorems for quadratic forms in
% random
% variables having long-range dependence. {\it Probab. Theory Related Fields}
% {\bf 74}, 213-240.
% \end{description}

 \bibliographystyle{elsarticle-harv}
\bibliography{biblio}

%%%%%%%%%%%%%%%%%%%%%%%%%%%%%%%%%%%%%%%%%%%%%%%%%%%%%%%%%%%%%%%%%%%%%%%%%%%%%%%%%%%%%%%%%%%%%%%%%%%%%%%%%%%%%%%%%%%%%%%%%%%%

% \vskip .65cm
% \noindent
% Istituto Dalle Molle di Studi sull'Intelligenza Artificiale (IDSIA) Manno, Switzerland
% \vskip 2pt
% \noindent
% E-mail: alessio@idsia.ch, francesca@idsia.ch, zaffalon@idsia.ch
% \vskip 2pt
% \noindent
% CNR IMATI, Milano, Italy
% \vskip 2pt
% \noindent
% E-mail: fabrizio@mi.imati.cnr.it
% \vskip .3cm

%%%%%%%%%%%%%%%%%%%%%%%%%%%%%%%%%%%%%%%%%%%%%%%%%%%%%%%%%%%%%%%%%%%%%%%%%%%%%%%%%%%%%%%%%%%%%%%%%%%%%%%%%%%%%%%%%%%%%%%%%%%%
%%%%%%%%%%%%%%%%%%%%%%%%%%%%%%%%%%%%%%%%%%%%%%%%%%%%%%%%%%%%%%%%%%%%%%%%%%%%%%%%%%%%%%%%%%%%%%%%%%%%%%%%%%%%%%%%%%%%%%%%%%%%

\end{document}